\documentclass[aos,seceqn,nameyear,dvips]{arximspdf}
\usepackage{graphics}
\doi{10.1214/09-AOS753}
\volume{38}
\issue{3}
\pubyear{2010}
\firstpage{1403}
\lastpage{1432}

\makeatletter

\newtheorem{Theorem}{Theorem}[section]
\newtheorem{Corollary}{Corollary}[section]

\newproclaim{Remark}{Remark}[section]

\makeatother

\begin{document}
\begin{frontmatter}

\title{Estimation in additive models with highly or nonhighly
correlated covariates}
\runtitle{Highly correlated additive modeling}

\begin{aug}
\author[A]{\fnms{Jiancheng} \snm{Jiang}\corref{}\thanksref{m1}\ead[label=e1]{jjiang1@uncc.edu}},
\author[B]{\fnms{Yingying} \snm{Fan}\thanksref{m2}\ead[label=e2]{fanyingy@marshall.usc.edu}} and
\author[C]{\fnms{Jianqing} \snm{Fan}\thanksref{m3}\ead[label=e3]{jqfan@princeton.edu}}
\runauthor{J. Jiang, Y. Fan and J. Fan}
\affiliation{University of North Carolina at Charlotte, University of
Southern California
and~Princeton University}
\address[A]{J. Jiang\\
Department of Mathematics and Statistics\\
University of North Carolina at Charlotte\\
Charlotte, North Carolina 28223\\
USA\\
\printead{e1}}
\address[B]{Y. Fan\\
Information and Operations\\
\quad Management Department\\
Marshall School of Business\\
University of Southern California\\
Los Angeles, California 90089\\
USA\\
\printead{e2}}
\address[C]{J. Fan\\
Department of Operations Research \\
\quad and Financial Engineering\\
Princeton University\\
Princeton, New Jersey 08544\\
USA\\
\printead{e3}}
\end{aug}

\thankstext{m1}{Supported by NSF Grant DMS-09-06482.}
\thankstext{m2}{Supported by NSF Grant DMS-09-06784.}
\thankstext{m3}{Supported by NSF Grants DMS-07-04337 and DMS-07-14554.}

\received{\smonth{1} \syear{2009}}
\revised{\smonth{8} \syear{2009}}

%
\begin{abstract}
Motivated by normalizing DNA microarray data and by predicting the
interest rates, we explore nonparametric estimation of additive models
with highly correlated covariates.
We introduce two novel approaches for
estimating the additive components, integration estimation and
pooled backfitting estimation. The former is designed for highly
correlated covariates, and the latter is useful for nonhighly
correlated covariates. Asymptotic normalities of the proposed
estimators are established.
Simulations are conducted to demonstrate
finite sample behaviors of the proposed estimators, and real data
examples are given to illustrate the value of the methodology.
\end{abstract}

%

\end{frontmatter}

\section{Introduction}\label{sec:intro}

The problem of estimating additive components with highly correlated
covariates arises from the normalization of DNA microarray. Since
the late 1980s, Affymetrix was founded with the revolutionary idea
to use semiconductor manufacturing techniques to create GeneChips
(an Affymetrix trademark) or generic DNA microarrays. It makes
quartz chips for the analysis of DNA microarrays and covers about
$82\%$ of the DNA microarray market.
A single chip can be used to do thousands of experiments in
parallel, so it produces a lot of Affymetrix GeneChip arrays which
demand proper normalization for removing systematic biases such as
the intensity effects.

Much research has been devoted to eliminating the systematic biases
such as the dye, intensity and print-tip block effects. Examples
include the rank-invariant selection method of Tseng et al. (\citeyear
{Tsengetal01}),
the lowess method of Dudoit et al. (\citeyear{Dudoitetal02}) and various
information aggregation methods of Fan et al. (\citeyear{Fan05a}),
Fan, Huang and Peng (\citeyear{FHP05}), Huang, Wang
and Zhang (\citeyear{HWZ05}) and Huang and Zhang (\citeyear{HZ05}),
among others.

Fan et al. (\citeyear{Fan05a}),
Fan, Huang and Peng (\citeyear{FHP05}) propose a semilinear in-slide
model (SLIM) to
remove intensity effects and identify significant genes for
Affymetrix arrays. Suppose that there are $G$ genes and for each
gene there are $J$ replications ($J\ge2$). Let $A_{gj}$ and
$B_{gj}$ be the log-detection signal of the $g$th probe set in
the $j$th control and treatment arrays, respectively.
Then, we compute the log
intensities and log-ratios, respectively, as
\[
X_{gj} = ({A}_{gj} + B_{gj} )/2,\qquad Y_{gj} = B_{gj} - {A}_{gj}.
\]
%
Fan et al. (\citeyear{Fan05a}),
Fan, Huang and Peng (\citeyear{FHP05}) use the following model to
estimate the
treatment effect, the smooth intensity effect:
%
%
\begin{equation} \label{b1}
Y_{gj}=\alpha_g+m_j(X_{gj})+\varepsilon_{gj},\qquad g=1,\ldots,G;
j=1,\ldots,J,
\end{equation}
where $\alpha_g$ is the treatment effect on gene $g$,
$m_j(X_{gj})$ represents the array-dependent intensity effect to be
estimated and $\varepsilon_{gj}$'s are independent noises with zero
means. For identifiability, we assume that $E[m_j(X_{gj})]=0$.

Directly estimating the treatment effects $\{\alpha_g\}$ is not a good idea
due to the existence of unknown intensity effects, as well as the small
size $J$.
In this paper we first treat $\{\alpha_g\}$ as nuisance parameters and
focus on the estimation of $m_j$'s.
Once a good estimate $\hat m_j$ of $m_j$ for each $j$ is obtained,
$\alpha_g$ can be estimated as
$\hat\alpha_g = \frac{1}{J}\sum_{j=1}^J(Y_{gj} - \hat m_j(X_{gj}))$.
Therefore, it is essential to efficiently estimate treatment effects $\{
\alpha_g\}$.
The setup applies to the c-DNA microarray data
[Fan, Huang and Peng (\citeyear{FHP05}), Huang and Zhang (\citeyear
{HZ05})] and Agilent microarray data
[Patterson et al. (\citeyear{Petal06})]. Moreover, it is also
applicable to other problems
where confounding effects can nonparametrically be removed.

\begin{figure}

\includegraphics{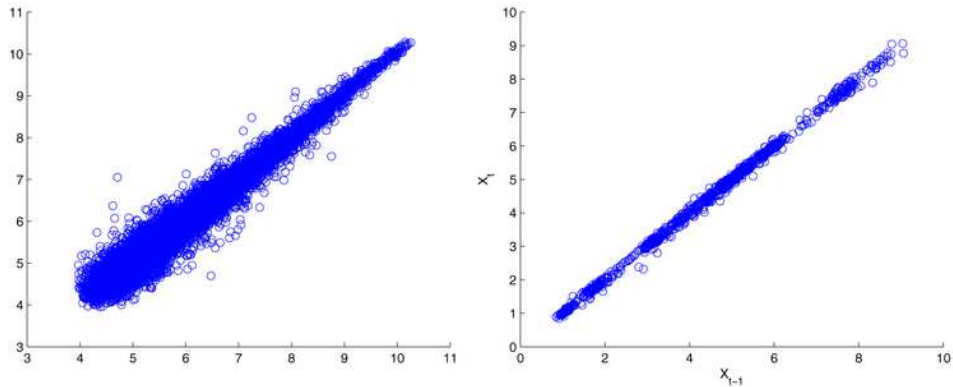}

\caption{Left panel: highly correlated log intensities in
Affymetrix array data with $J=2$; right panel: highly correlated
interest rates. $\{X_t\}$ represents the weekly data of the 6-month
treasury bill secondary market rates in the period of June 1, 1988 to
June 1, 2008.}\label{fig0}
\end{figure}

Fan et al. (\citeyear{Fan05a}) used a backfitting algorithm to estimate
iteratively the intensity effect and the treatment effect. While
this method is successful for removing the systematic biases in some
certain situations, mathematical properties of the resulting
estimators are unknown which requires further study of the
estimation. On the other hand, when performing the estimation
method, we found that it is unstable and even fails to converge in
some situations. A careful study of this problem reveals that it is
caused by the high correlation between intensities.
An illustrating example is the DNA microarrays data analyzed in
Fan et al. (\citeyear{Fan05a}). In this example, the log-intensities across
different chips are highly correlated, which is evidenced in Figure
\ref{fig0}(left), due to the repeatability and accuracy of the
measurements. A close look at the almost identical relationship
between covariates suggests that $|X_{g1}-X_{g2}|\to0$. This
suggests a simple working model for calibrating the following plausible
correlation structure:
\[
X_{g1}=X_{g2}+b_Gu_{g2},
\]
where $b_G\to0$ and $u_{g2}$ is random noise.
Under such a setting, the problem of effectively estimating the confounding
effect $\{m_j(\cdot)\}$ challenges statisticians. The high correlation
reduces the accuracy of estimating $m_j(\cdot)$, but $G$ in such an
application is also very large, in an order of tens of thousands.

%

The problem of highly correlated covariates appears often in
modeling time series data such as interest rates. Suppose that we
would like to use the past 4 weeks' ($X_{t-1}, \ldots, X_{t-4}$) interest
rates to forecast the return of a stock or index $Y_t$ or the
interest rate itself $Y_t = X_t$ in the next week. A reasonable
nonparametric model is the following additive model:
\[
Y_t = \mu+ m_1(X_{t-1}) + \cdots+ m_4 (X_{t-4}) + \varepsilon_t.
\]
Due to the continuity of the interest rate dynamics, the covariates
in the above additive model is also highly correlated and can be
handled by the idea in this paper. Figure \ref{fig0}(right) shows the scatter
plot of $X_t$ versus $X_{t-1}$ using the weekly data for the 6-month
treasury bill secondary market rates in the period of June 1, 1988 to
June 1, 2008.

Existing methods in the literature do not appear enough to address the
problem with additive modeling with highly correlated covariates,
and a new methodology is needed. In particular, in addition to the
aforementioned failure in convergence, the backfitting algorithm
usually converges slowly due to the very
large number of genes $G$ which is usually in the order of tens of
thousands in a typical microarray application. This motivates us to
develop statistical methods fitting the smooth confounding effect
model (\ref{b1}) with/without highly correlated intensity effects.

The above model received attention in Fan et al. (\citeyear{Fan05a}),
Fan, Huang and Peng (\citeyear{FHP05}). However, there is no formal
study of modeling highly correlated covariates~$X_{gj}$. For the usual
correlation situation, Fan, Huang and Peng (\citeyear{FHP05})
considered the estimator of $m_j$ using the profile least squares and
obtained only an upper bound for the conditional mean squared error of
the estimator. However, information across arrays is not used, and the
asymptotic distribution of the estimator is unknown which makes the
inference about the intensity effects difficult.

In this investigation, we introduce two methods for estimating the
nonparametric components $m_j$, integration estimation and pooled
backfitting estimation. The former is tailored for modeling highly
correlated intensity effects and is a
noniterative estimator with fast implementation. It relies on
estimating the derivative
function in a varying coefficient model, and allows us to handle a
very large amount of observations.
The latter is an iterative estimate which is designed for modeling
nonhighly correlated intensity effects.
Asymptotic normalities of
the proposed estimators are established. The extent to which the high
correlation affects the rates of convergence is explicitly given.
Simulation studies are conducted to demonstrate finite sample behaviors
of the proposed
methods. 

The paper is organized as follows.
In Section \ref{sec:c} we introduce the integration estimation method
along with an alternative of robustness.
In Section \ref{sec:b} we develop
pooled backfitting estimation of the intensity effects.
In Section \ref{sec:d} we conduct simulations.
In Section \ref{sec:e} we illustrate the proposed methodology by two
real data examples.
Finally we conclude the
paper with a discussion. Details of assumptions and proofs of
theorems are given in Appendices \ref{appI} and \ref{appII}.

\section{Estimation of additive components when covariates are highly
correlated}\label{sec:c}

To use information across arrays, one can take a difference operator
to remove the nuisance parameters $\{\alpha_g\}$ which leads to
additive models. Specifically, let
$Y_g^{(k)}=Y_{g1}-Y_{gk}$ and
$\varepsilon_g^{(k)}=\varepsilon_{g1}-\varepsilon_{gk}$. Then by
(\ref{b1}), for $k=2,\ldots,J$,
%
%
\begin{equation}\label{b2}
Y_g^{(k)}=m_1(X_{g1})-m_k(X_{gk})+\varepsilon_{g}^{(k)},\qquad
g=1,\ldots,G,
\end{equation}
which are additive models introduced by Friedman and Stuetzle
(\citeyear{FS81}) and Hastie and Tibshirani (\citeyear{HT90}) where
$\varepsilon_{g}^{(k)}$ are the errors with zero means, and for $j\ne
k$,
$\operatorname{Cov}(\varepsilon_{g}^{(j)},\varepsilon_{g}^{(k)})=\sigma
^2$ and $\operatorname{Var}(\varepsilon_{g}^{(j)})=2\sigma^2$. The
additive components can be estimated via the backfitting method. Due to
the high correlation between $X_{g1}$ and~$X_{gk}$, the estimate based
on the backfitting algorithm usually fails in convergence, and the
existence of a backfitting estimator is problematic. Moreover,
asymptotic properties of the backfitting estimators are unknown in this
situation. Thus a new methodology is needed to deal with this problem.
To this end, in the following we focus on the cases with highly
correlated covariates and introduce the integration estimation and then
establish asymptotic normality of the resulting estimators under a
working model. The estimators are consistent, regardless of the working
model.

\subsection{Estimation when covariates are highly correlated}\label{sec:2.1}

As illustrated in the previous section, covariates $X_{gj}$ (for
$j=1,\ldots,J$) may be
very close and highly correlated, so it is convenient to assume that
%
%
\begin{equation} \label{con1}
\Delta_{gk}\equiv X_{g1}-X_{gk} \to0.
\end{equation}
Under such a setting, the asymptotic properties of the backfitting
estimates are unknown, and the convergence of the backfitting
algorithm may also be a problem since the required condition, that is,
the existence of the joint density of covariates, is not always
satisfied. See, for example, Opsomer and Ruppert (\citeyear{OR97},
\citeyear{OR98}). Assume
$m_1''$ is continuous; then by Taylor's expansion,
%
%
\begin{eqnarray} \label{cc1}
m_1(X_{g1})
&=& m_1(X_{gk})+m_1'(X_{gk})\Delta_{gk}\nonumber\\[-8pt]\\[-8pt]
&&{} + \tfrac{1}{2}m_1''(X_{gk})\Delta_{gk}^2+o(\Delta_{gk})^2.\nonumber
\end{eqnarray}
Substituting (\ref{cc1}) into (\ref{b2}), we obtain that
%
%
\begin{equation}\label{b2a}
Y_g^{(k)}
=m_{k1}(X_{gk})+m_1'(X_{gk})\Delta_{gk}+\tilde{\varepsilon
}_g^{(k)},
\end{equation}
where $m_{k1}(X_{gk}) = m_1(X_{gk}) -m_k(X_{gk})$ and
$\tilde{\varepsilon}_g^{(k)}=\frac{1}{2}m_1''(X_{gk})\Delta
_{gk}^2+o(\Delta_{gk})^2+{\varepsilon}_g^{(k)}$.
Model (\ref{b2a}) is actually a varying coefficient model, since the
coefficient functions $m_{k1}(\cdot)$ and $m'_1(\cdot)$ are unknown
functions of $X_{gk}$. This allows us to estimate the unknown
coefficient functions $m'_1(\cdot)$ using local smoothing
techniques.
Given an interior point $x\in\operatorname{supp}[f_k(\cdot)]$, using
the local
linear approximation when $|X_{gk}-x|\le h$,
we obtain that
%
%
\begin{eqnarray}\label{b6}
&&m_{k1}(X_{gk})+m_1'(X_{gk})\Delta_{gk}\nonumber\\[-8pt]\\[-8pt]
&&\qquad\approx
\alpha_0+\alpha_1(X_{gk}-x)+\Delta_{gk}\{\beta_0+\beta_1(X_{gk}-x)\}
.\nonumber
\end{eqnarray}
Then the coefficient function $m_1'(\cdot)$ can be estimated by minimizing
%
%
\begin{eqnarray}\label{b7a}
&&\sum_{g=1}^G\bigl[Y_g^{(k)}-\alpha_0-\alpha_1(X_{gk}-x)\nonumber\\[-8pt]\\[-8pt]
&&\hspace*{17.8pt}{} -\Delta_{gk}\{\beta_0+\beta_1(X_{gk}-x)\} \bigr]^2
K_h(X_{gk}-x),\nonumber
\end{eqnarray}
where $K_h(\cdot)=h^{-1}K(\cdot/h)$ with $K(\cdot)$ being a kernel
function and $h$ being a bandwidth
controlling the amount of data in smoothing.
Denote by
$\{\hat{\alpha}(x),\hat{\beta}(x)\}$ with $\hat{\alpha}(x) = (\hat
{\alpha}_0(x), \hat{\alpha}_1(x))$ and $\hat{\beta}(x) = (\hat{\beta
}_0(x), \hat{\beta}_1(x))$ the solution to the above
equation.
Then
$\hat{\beta}_0(x)$
and $\hat{\beta}_1(x)$ estimate $m_1'(x)$ and $m_1''(x)$,
respectively.
If $\Delta_{gk}=o(1)$, then
$E[\tilde{\varepsilon}_g^{(k)}|X_{gk}=x]=o(1)$, and hence
the above estimator is consistent. The method is noniterative and
can handle the situation where $G$ is very large.
Once the derivative $m_1'(\cdot)$ is given, the component $m_1$ in
model (\ref{b2}) can be derived
as follows.


Let $\mathbf K=\operatorname{diag}\{K_h(X_{1k}-x),\ldots,K_h(X_{Gk}-x)\}$,
$\hat{\theta}(x)=(\hat{\alpha}_0,\hat{\alpha}_1,\hat{\beta}_0,\hat
{\beta}_1)^T$,
\[
Z_g=\bigl(1,X_{gk}-x,b_Gu_{gk},b_Gu_{gk}(X_{gk}-x)\bigr)^T
\]
and
$\mathbf{Z}=(Z_1,\ldots,Z_G)^T$. Then $\hat{\theta}(x)$ admits the
following closed form:
%
%
\begin{eqnarray} \label{b7aa}\hat{\theta}(x)
=(\mathbf{Z}^T\mathbf K\mathbf{Z})^{-1}\mathbf{Z}^T\mathbf
K\mathbf{Y}_{(k)},
\end{eqnarray}
where
$\mathbf{Y}_{(k)}=(Y_1^{(k)},\ldots,Y_G^{(k)})^T$. Let
\[
\theta(x)=(m_{1k}(x),m'_{1k}(x),m_1'(x),m_1''(x))^T.
\]
Then
$\hat{\theta}(x)$ estimates $\theta(x)$, and $m_1'(x)$ is estimated
by $\hat{m}_1'(x;k)=e_3^T\hat{\theta}(x)$
with $e_3=(0,0,1,0)^T$.

Since averaging can reduce the variance of estimation,
we propose to estimate $m_1'(\cdot)$ by the following average:
%
%
\begin{equation} \label{b10}
\hat{m}_1'(x)=(J-1)^{-1}\sum_{k=2}^J\hat{m}_1'(x;k).
\end{equation}
Note that, for each $k$, $\hat{m}_1'(x;k)$ is consistent. The estimator
$\hat{m}_1'(\cdot)$
is also consistent.
From the estimated derivative function, the original function
$m_1(\cdot)$ can
consistently be estimated using integration which we now detail below.

Let $F_j(\cdot)$ and $f_j(\cdot)$ be, respectively, the distribution
and density functions of~$X_{gj}$.
Due to the identifiability condition $E[m_1(X_{g1})]=0$ and
$
m_1(x)=m_1(x_0)+\int_{x_0}^xm_1'(t) \, dt (\mbox{for
any } x_0\in\operatorname{supp}[F_1(\cdot)]) $, we obtain that
\[
\int
\biggl\{m_1(x_0)+\int_{x_0}^xm_1'(t) \,dt \biggr\} \,dF_1(x)=0
\]
and hence
$
m_1(x_0)=-\int\int_{x_0}^xm_1'(t) \,dt \,dF_1(x).
$ 
Therefore, $m_1(x)$ can be estimated by
%
%
\begin{equation} \label{b8}
\hat{m}_1(x)=-\int\int_{x_0}^x\hat{m}_1'(t) \,dt \,d\hat{F}_1(x)
+\int_{x_0}^x\hat{m}_1'(t) \,dt,
\end{equation}
where $\hat{F}_1$
is the empirical estimator of $F_1$. Note that the first term in
(\ref{b8}) is a constant, making merely the estimated function to
satisfy an empirical version of the identifiability condition.
Similarly, we can estimate the other components' $m_j$'s (for
$j=2,\ldots,J$) in model (\ref{b2}). Such defined estimators are
naturally consistent due to consistency of
the estimators of derivative functions.

\subsection{Asymptotic normality}

To provide in-depth analysis on the behavior of the estimators
defined in (\ref{b7aa})--(\ref{b8}), we model explicitly the
high correlation among covariates. One viable choice is to employ
the following working model:
%
%
\begin{equation} \label{b3}
X_{g1}=X_{gk}+b_Gu_{gk},
\end{equation}
where $b_G\rightarrow0$ and $\{u_{gk}\}_{g=1}^G$ are noises of zero
mean and finite variance.
Assume that the density
function of $u_{gk}$, $p_k(x)$, has a compact support and that $\{
u_{gk}\}$ are independent of $\{X_{gk}\}$
for fixed $k$.
This specification allows for heteroscedasticity of the errors.
Obviously, in
model (\ref{b3}) the correlation between $X_{g1}$ and $X_{gk}$ goes
to one as $b_G\to0$. There are various alternative methods for
modelling high correlation between two variables.
We focus only on model (\ref{b3}) to make an attempt.
Note that the working model (\ref{b3}) is only used to derive the
asymptotic properties. The estimator itself does not depend on such
an assumption.

Denote by
$\mu_j(K)=\int t^jK(t) \,dt$ and $\nu_j(K)=\int t^j K^2(t) \,dt$.
Let
$\mathbf{H}=\break\operatorname{diag}(\mathbf h,b_G\mathbf h)$,
$\mathbf{S}=\operatorname{diag}(\mathbf{N},\mathbf{N})$,
$\mathbf{V}=\operatorname{diag}(\bolds\nu,\bolds\nu)$,
$\mathbf{C}=\operatorname{diag}(\mathbf{c}_2,\mathbf{c}_2)$
and
$\mathbf{c}^*=(\mathbf{c}_0^T,\mathbf{c}_0^TE(u_{1k}^3))^T$
where
$\mathbf h=\operatorname{diag}(1,h)$,
$\mathbf{N}=\operatorname{diag}\{\mu_0(K),\mu_2(K)\}$,
$\bolds\nu=\break\operatorname{diag}\{\nu_0(K),\nu_2(K)\}$
and
$\mathbf{c}_j=(\mu_j(K),\mu_{j+1}(K))^T$.
The following theorems describe the asymptotic properties of the
proposed estimators.
\begin{Theorem}\label{Th1}
Suppose that the conditions in Appendix \ref{appI} hold. Under the working model
(\ref{b3}),
if $Gh^5=O(1)$ and $Ghb_G^4=O(1)$, then
\[
\sqrt{Gh} \bigl\{ \mathbf{H}[\hat{\theta}(x)-\theta(x)]-\mathbf
{b}(x)\bigl(1+o_p(1)\bigr) \bigr\}
\stackrel{\mathcal D}{\longrightarrow} {\mathcal N}(0,\Sigma(x)),
\]
where
\[
\mathbf{b}(x)=\tfrac{1}{2}h^2\mathbf{S}^{-1}\mathbf
{C}\bigl(m_{1k}''(x),b_Gm_1^{(3)}(x)\bigr)^T
+\tfrac{1}{2}b_G^2m_1''(x)\mathbf{S}^{-1}\mathbf{c}^*
\]
and
$\Sigma(x)=2\sigma^2f_1^{-1}(x)\mathbf{S}^{-1}\mathbf{V}\mathbf{S}^{-1}$.
\end{Theorem}
\begin{Corollary}\label{col1}
Under the conditions in Theorem \ref{Th1},
\[
\sqrt{Gh}b_G
\{\hat{m}_1'(x;k)-m_1'(x) -b_1(x) \} \stackrel{\mathcal
D}{\longrightarrow} {\mathcal N}(0,\sigma_1^2(x)),
\]
where
\[
b_1(x)=\tfrac{1}{2}h^2\mu_2(K)\mu_0^{-1}(K)m_1^{(3)}(x)\bigl(1+o_p(1)\bigr)
+\tfrac{1}{2}b_G m_1''(x)E(u_{1k}^3)\bigl(1+o_p(1)\bigr)
\]
and
$\sigma_1^2(x)=2\sigma^2f_1^{-1}(x)\nu_0(K)\mu_0^{-2}(K)$.
\end{Corollary}

The above corollary shows that the data from two arrays suffice to
obtain a consistent estimate of the derivative function. However,
the high correlation reduces the effective sample size from $G$
to $Gb_G^2$, in terms of the rates of convergence.

In order to present asymptotics of the average estimator (\ref{b10}),
we need the dependence structure of $\{u_{gk}\}$ across $k$.
Let
$\rho(\ell,k)=E(u_{g\ell}u_{g,k})$, which does not depend on $g$,
and
\[
\rho=\Biggl\{\sum_{k=2}^J\rho(k,k)
+\sum_{k_1=2}^{J}\sum_{k_2=2}^J\rho(k_1,k_2)\Biggr\}\bigg/(J-1)^2.
\]
\begin{Theorem}\label{Th1a}
Under the conditions in Theorem \ref{Th1},
\[
\sqrt{Gh}b_G
\{\hat{m}_1'(x)-m_1'(x) -b_1(x) \} \stackrel{\mathcal
D}{\longrightarrow} {\mathcal N}(0,\sigma_2^2(x)),
\]
where
$\sigma_2^2(x)= \rho\sigma^2f_1^{-1}(x)\nu_0(K)\mu_0^{-2}(K)$.
\end{Theorem}

The above asymptotics of the estimators is derived under the working
model (\ref{b3}).
However, as previously stated, if condition (\ref{con1}) holds,
our estimator for $m_1'(x)$ is consistent
whether or not the working model (\ref{b3}) holds.
This furnishes robustness of our estimator $\hat{m}'(x)$ against
mis-specification of the correlation between covariates.
If interested in estimating the derivative function,
one can directly compute the asymptotic bias and variance of
$\hat{m}_1'(x)$ and obtain the optimal bandwidth by minimizing the
asymptotic mean square error so that a data-driven
bandwidth selection rule can be developed as in the one-dimensional
nonparametric regression problem.
In the following we state the
asymptotic normality of the integrated estimator.
\begin{Theorem}\label{Th2}
Suppose that the conditions in Appendix \ref{appI} hold. Under the working model
(\ref{b3}),
if $Gb_G^2h^4=O(1)$ and $Gb_G^4=O(1)$, then
%
\[
\sqrt{G}b_G\bigl\{\hat{m}_1(x)-m_1(x)
-B_1(x)\bigl(1+o_p(1)\bigr)\bigr\}
\stackrel{{\mathcal D}}{\longrightarrow}{\mathcal
N}(0,\sigma^2(x)),
\]
where
$\sigma^2(x)=\frac{1}{4}\rho\sigma^2f_1^{-2}(x)$
and
\begin{eqnarray*}
B_1(x)&=&\tfrac{1}{2}h^2\mu_2(K)\mu_0^{-1}(K)\{
m_{1}''(x)-E[m_1''(X_{11})]\}\\
&&{}+\tfrac{1}{2}b_G E(u_{1k}^3)\{m'_1(x)-E[m_1'(X_{11})]\}.
\end{eqnarray*}
\end{Theorem}
\begin{Remark}
If $h=o(\sqrt{b_G})$,
then the bias term is
\[
B_1(x)=\tfrac{1}{2}b_G E(u_{1k}^3)\{m'_1(x)-E[m_1'(X_{11})]\}\bigl(1+o(1)\bigr)
\]
and hence the asymptotic normality of the estimator does not depend on
the smoothing parameter $h$ nor the kernel $K$.
It parallels the result of Jiang, Cheng and Wu (\citeyear{JCW02}) for estimating
distribution functions
and contrasts with the dependence on smoothing parameter of the
nonparametric function estimation.
\end{Remark}
\begin{Remark}
The estimate $\hat{m}_1(\cdot)$ achieves a maximum convergence rate
$O(G^{-1/4})$ when $b_G=O(G^{-1/4})$.
The convergence rate can be improved if one uses a higher order
polynomial approximation in (\ref{b6}).
\end{Remark}


\subsection{A pooled robust approach}\label{sec:s}

In model (\ref{b2}), we aim at estimating $m_1(\cdot)$. It has
various versions of implementations. To illustrate the idea, we use
aggregated local constant approximation along with the $L_1$-loss to
illustrate the versatility. For $|X_{gk}-x|=O(h)$, we have
$m_{k1}(X_{gk})\approx m_{k1}(x)$ and $m_1'(X_{gk})\approx m_1'(x)$.
Then, by~(\ref{b2a}), we can run the local regression by minimizing
%
%
\begin{equation}\label{eq1a}
\sum_{k=2}^J\sum_{g=1}^G \bigl|Y_g^{(k)}-\alpha_{k,0}
-\beta_0\Delta_{gk} \bigr| K_h(X_{gk}-x)
\end{equation}
with $Y_g^{(k)}=Y_{g1}-Y_{gk}$, $\Delta_{gk} = X_{g1} - X_{gk}$ and
$K_h(\cdot)=h^{-1}K(\cdot/h)$. Notice that we pool data from
different
replicates in (\ref{eq1a}) to obtain more accurate estimators,
and the $L_1$ norm is used to alleviate the influence of outliers.
Denote by $(\hat{\alpha}_{2,0}, \ldots,\hat{\alpha}_{J,0}, \hat{\beta
}_0(x))$ the solution to the
above minimization problem. Then $\hat{\beta}_0(x)$ estimates
$m_1'(x)$. Integrating $\hat{\beta}_0(x)$ leads to an estimate of
$m_1(x)$. In our experience, this estimation approach performs
similarly to the method in previous sections.

\section{Backfitting estimation of additive components}\label{sec:b}

In this section, we introduce pooled backfitting estimators of $m_j$
and study their asymptotic properties under nonhigh correlation
situations.

\subsection{Fitting a bivariate additive model using the local linear
smoother based on the backfitting algorithm}

There are some methods for fitting the additive model~({\ref{b2}).
For example, the common backfitting estimation of Buja, Hastie
and Tibshirani (\citeyear{BHT89}) and Opsomer and Ruppert (\citeyear
{OR97}, \citeyear{OR98}),
the marginal integration methods of
Tj\o theim and Auestad (\citeyear{TA94}), Linton and Nielsen (\citeyear
{LN95}) and
Fan, H\"{a}rdle and Mammen (\citeyear{FHM98}),
the estimating equation method of Mammen, Linton and Nielsen (\citeyear{MLN99})
and the smooth backfitting method in Nielsen and Sperlich
(\citeyear{NS05}), among others. For illustration, we will use the common
backfitting algorithm based on the local linear smoother as a building
block to estimate the additive components. Other estimation methods
can similarly be applied.

To ensure identifiability of the additive component functions
$m_j(\cdot)$, we impose the constraint
$E[m_j(X_{gj})]=0$ for $j=1,\ldots,J$.
Fitting the additive component $m_j(\cdot)$ in (\ref{b2}) requires
choosing bandwidths $\{h_j\}$. The optimal choice of $h_j$ can be
obtained as in Opsomer and Ruppert (\citeyear{OR98}).
We
here follow notation that was introduced by Opsomer and Ruppert
(\citeyear{OR97}). Put
$K_{h_j}(x)=h_j^{-1}K(\frac{x}{h_j})$, $K_s(v)=v^{s-1}K(v)$,
${\mathbf H}_j=\operatorname{diag}(1,h_j)$, ${\mathbf
m}_j=\{m_j(X_{1j}),\ldots,m_j(X_{Gj})\}^T$,
$\mathbf{X}_j=(X_{1j},\ldots,X_{Gj})^T$
and
$\mathbf{Y}_k=(Y_{1}^{(k)},\ldots,Y_{G}^{(k)})^T$.
The smoothing matrices for local polynomial
regression are
\[
{\mathbf{S}}_j=({\mathbf s}_{j,X_{1j}},\ldots,{\mathbf
s}_{j,X_{Gj}})^T,
\]
where ${\mathbf s}_{j,x_j}^T$ represents the equivalent kernel
for the $j$th covariate at the point $x_j$.
%
%
\begin{equation} \label{B1}
{\mathbf s}_{j,x_j}^T ={\mathbf e}_1^T\bigl({{\mathbf X}_{x_j}^{(j)}}^T
{\mathbf K}_{x_j} {\mathbf X}_{x_j}^{(j)}\bigr)^{-1} {{\mathbf
X}_{x_j}^{(j)}}^T {\mathbf K}_{x_j}.
\end{equation}
Here $\mathbf e_1^T=(1,0)$,
${\mathbf K}_{x_j}=\operatorname{diag}\{K_{h_j}(X_{1j}-x_j),\ldots
,K_{h_j}(X_{Gj}-x_j)\}$
and
\[
{\mathbf X}_{x_j}^{(j)}=\left[
\matrix{
1 & (X_{1j}-x_j)\cr
\vdots& \vdots\cr
1 & (X_{Gj}-x_j)}
\right].
\]

From (\ref{b2}), ${\mathbf m}_j$'s can be estimated through the
solutions to the
following set of normal equations [see Buja, Hastie and Tibshirani
(\citeyear{BHT89}), Opsomer and Ruppert (\citeyear{OR97})]:
\[
\left[\matrix{
{\mathbf I}_G & \mathbf{S}_1^*\cr
\mathbf{S}_k^* & {\mathbf I}_G
}\right]
\left[\matrix{
\hat{\mathbf m}_1\cr
-\hat{\mathbf m}_k}
\right] = \left[\matrix{
{\mathbf{S}}_1^*\cr
{\mathbf{S}}_k^*}
\right] \mathbf{Y}^{(k)},
\]
where ${\mathbf{S}}_j^*=({\mathbf I}_G-{\mathbf1}{\mathbf1}^T/G){\mathbf{S}}_j$
is the centered smoother matrix, and ${\mathbf1}$ is a $G\times1$
vector whose
elements are all ones. In practice, the backfitting
algorithm [Buja, Hastie and Tibshirani (\citeyear{BHT89})] is usually
used to solve these
equations, and the backfitting estimators converge to the solution,
%
%
\begin{equation} \label{B2a}
\left[\matrix{
\hat{\mathbf m}_1^{(k)}\cr
-\hat{{\mathbf m}}_k}
\right] = \left[\matrix{
{\mathbf I}_G & {\mathbf{S}}_1^* \cr
{\mathbf{S}}_k^* & {\mathbf I}_G
} \right]^{-1}
\left[\matrix{
{\mathbf{S}}_1^*\cr
{\mathbf{S}}_k^*}
\right] \mathbf{Y}^{(k)}\qquad\mbox{for } k=2,\ldots,J,
\end{equation}
where the superscript in $\hat{\mathbf m}_1^{(k)}$ is used to
stress the
dependence of $\hat{\mathbf m}_1$ on $k$.

If $\Vert\mathbf{S}_1^* \mathbf{S}_k^*\Vert<1$, then the
backfitting estimators exist and are unique where
we use
$\|\mathbf{A}\|$ to denote the maximum row sum matrix norm of the square
matrix $\mathbf{A}\dvtx
\|\mathbf{A}\|=\max_{1\le i\le G}\sum_{j=1}^G|A_{ij}|$. A sufficient condition
for $\Vert\mathbf{S}_1^* \mathbf{S}_k^*\Vert<1$ is
%
%
\begin{equation}\label{2.3a}
\sup_{x_1,x_k} \biggl|\frac{f_{1k}(x_1,x_k)}
{f_{1}(x_1)f_{k}(x_k)}-1 \biggr|<1,
\end{equation}
where $f_j(x_j)$ is the density of $X_j$, and $f_{1k}(x_1,x_k)$ is
the joint density of $X_{g1}$ and $X_{gk}$ [see Opsomer and Ruppert
(\citeyear{OR97})]. We assume in this
section the above condition holds. Note that this condition does not
hold for the working model (\ref{b3}) since the joint density of
$(X_{g1},X_{gk})$ is nearly degenerate. Solving (\ref{B2a}),
we get
%
%
\begin{equation}\label{B3}
\hat{{\mathbf m}}_1^{(k)}=\{{\mathbf I}_G-({\mathbf
I}_G-\mathbf{S}_1^*\mathbf{S}_k^*)^{-1}({\mathbf I}_G-\mathbf{S}_1^*)\}
\mathbf{Y}^{(k)}\equiv
\mathbf{W}_{1k}\mathbf{Y}^{(k)}.
\end{equation}

Since averaging can reduce the variance, we propose to estimate
$\mathbf m_1$ by
%
%
\begin{equation} \label{B5}
\hat{\mathbf m}_1=(J-1)^{-1}\sum_{k=2}^J
\hat{\mathbf m}_1^{(k)},
\end{equation}
which is termed as the pooled backfitting estimator of~$\mathbf m_1$.
For other components~$\mathbf m_j$,
they can be estimated in a similar way. Thus, in the following, we
will focus on the estimation of $\mathbf m_1$. The integration method in
the previous section is simpler and much faster to compute since it
uses only one
smoothing parameter $h$ and does not
involve any iteration.

To derive the asymptotic properties of $\hat{\mathbf m}_1$, in the
following we introduce some notation in Opsomer and Ruppert (\citeyear{OR97}).
Define
\[
D_{x,h_1}=\{t\dvtx(x+h_1t)\in\operatorname{supp}(f_1)\}\cap
\operatorname{supp}(K).
\]
Then $x$ is called ``an interior point'' if any only if
$D_{x,h_1}=\operatorname{supp}(K)$. Otherwise, $x$ is a boundary point.
Define the kernel $K_{(1)}(u) =
K(u)/\mu_0(K)$, which is the asymptotic
counterpart of the equivalent kernel induced by the local linear
fit.
Then $\mu_2(K_{(1)})=\mu_2(K)\mu_0^{-1}(K)$
and
$\nu_0(K_{(1)})=\nu_0(K)\mu_0^{-2}(K)$.
Let $T_{1k}^*$ be a matrix whose $(i,j)$th element
is
\[
[T_{1k}^*]_{ij}=G^{-1}\{f_{1k}(x_1,x_k)f_1^{-1}(x_1)f_k^{-1}(x_k)-1\}.
\]
Let $\mathbf t_g^T$ represent the $g$th row of $(I-T_{1k}^*)^{-1}$, and
$\mathbf e_g$ be the $g$th unit vector.
\begin{Theorem}\label{thb}
Suppose that the conditions in Appendix \ref{appI} hold. If $X_{g1}$ is an
interior point,
then as $G\to\infty$:
\begin{longlist}
\item the bias of $\hat{m}_1(X_{g1})$ conditional on $\mathbf
{X}=(\mathbf{X}
_1,\mathbf{X}_k)$ is
\[
E\{\hat{m}_1(X_{g1})-m_1(X_{g1}) | \mathbf{X}\}
=b_1-b_2+O_p\biggl(\frac{1}{\sqrt{G}}\biggr)+o_p\Biggl(\sum_{j=1}^Jh_j^2\Biggr),
\]
where
$b_1=\frac{1}{2}h_1^2\mu_2(K_{(1)})[m_1''(X_{g1})+\{(\mathbf
t_g^T-\mathbf e
_g^T)\mathbf m_1''-E(m_1''(X_{g1}))\}]$
and
\[
b_2=\frac{1}{2}\mu_2\bigl(K_{(1)}\bigr)\frac{1}{J-1}\sum_{k=2}^Jh_k^2 \{\mathbf t
_g^TE(m_k''(X_{gk}) | \mathbf{X}_1)-E(m_k''(X_{gk})) \};
\]

\item the variance of $\hat{m}_1(X_{g1})$ conditional on $\mathbf
{X}$ is
\[
\operatorname{Var}\{\hat{m}_1(X_{g1}) |\mathbf{X}\}=\frac{J}{J-1}\frac
{1}{Gh_1}\sigma
^2f_1^{-1}(X_{g1})\nu_0\bigl(K_{(1)}\bigr)+o_p\biggl(\frac{1}{Gh_1}\biggr).
\]
\end{longlist}
\end{Theorem}

As in Corollary 4.3 of Opsomer and Ruppert (\citeyear{OR97}),
if the covariates are independent, the conditional bias of $\hat
{m}_1(X_{g1})$ in the interior of $\operatorname{supp}(f)$ can be approximated
by
\begin{eqnarray*}
E\{\hat{m}_1(X_{g1})-m_1(X_{g1}) |\mathbf{X}\}
&=&\frac{1}{2}h_1^2\mu_2\bigl(K_{(1)}\bigr)\{m_1''(X_{g1})-E(m_1''(X_{g1}))\}\\
&&{}+O_p\bigl(1/\sqrt{G}\bigr)+o_p\Biggl(\sum_{j=1}^Jh_j^2\Biggr).
\end{eqnarray*}

\subsection{Fitting a J-variate additive model using local linear
smoother based on backfitting}

In the previous section, we used the differences between any two
different replicates for genes
to eliminate the nuisance parameters.
It resulted in two-dimensional additive models, which were easy to
implement, but for each additive model, the estimator was asymmetric.
In the following we use differences between any replicate and the
average of those replicates.
This will lead to a \mbox{$J$-dimensional} additive model with symmetric estimation.

Let
\[
\bar{Y}_g=J^{-1}\sum_{j=1}^J Y_{gj},\qquad
\bar{m}(X_g)=J^{-1}\sum_{j=1}^J m_{j}(X_{gj})
\]
and
$\bar{\varepsilon}_g=J^{-1}\sum_{j=1}^J \varepsilon_{gj}$. Then by
(\ref{b1}) we have
%
%
\begin{equation}\label{b2a1}
\bar{Y}_g=\alpha_g+\bar{m}(X_g)+\bar{\varepsilon}_g.
\end{equation}
Subtracting (\ref{b2a1}) from (\ref{b1}), we obtain that for
$j=1,\ldots,J$,
%
%
\begin{equation}\label{b2b}
Y_{gj}^*=-\frac{1}{J}\sum_{k\ne
j}m_k(X_{gk})+\frac{J-1}{J}m_j(X_{gj})+\varepsilon_{gj}^*,
\end{equation}
where $Y_{gj}^*=Y_{gj}-\bar{Y}_g$ and
$\varepsilon_{gj}^*=\varepsilon_{gj}-\bar{\varepsilon}_g$. It can be
seen that $\operatorname{Var}(\varepsilon^*_{gj})=(1 - 1/J) \sigma^2$ and
$\operatorname{Cov}(\varepsilon^*_{gj},\varepsilon^*_{kj})=0$ for $g\ne
k$. For
any fixed $j$, let
\[
m_{j,j}^*(X_{gj})=(J-1)J^{-1}m_j(X_{gj})
\]
and
$m_{k,j}^*(X_{gk})=-J^{-1}m_k(X_{gk})$ for $k\ne j$. Then
(\ref{b2b}) becomes
%
%
\begin{equation}\label{b2b1}
Y_{gj}^*=\sum_{k=1}^Jm_{k,j}^*(X_{gk})+\varepsilon_{gj}^*.
\end{equation}
This is a $J$-variate additive model. Again, we can estimate
the additive components using the local linear smoother based on the
backfitting algorithm.

Fitting the additive component $m_j$ in (\ref{b2b}) requires
choosing bandwidths $\{h_j\}$. The optimal choice of $h_j$ can be
obtained as in Opsomer and Ruppert (\citeyear{OR98}) and Opsomer
(\citeyear{O00}).
Put
\[
{\mathbf m}_{k,j}^*=\{m_{k,j}^*(X_{1k}),
\ldots,m_{k,j}^*(X_{Gk})\}^T
\quad\mbox{and}\quad
{\mathbf
Y}_j^*=(Y_{1j}^*,\ldots,Y_{Gj}^*)^T.
\]
Then the additive components
can be estimated through the solutions to the following set of
normal equations:
\[
\left[\matrix{
{\mathbf I}_G & {\mathbf S}_1^* & \cdots& {\mathbf S}_1^*\cr
{\mathbf S}_2^* & {\mathbf I}_G & \cdots& {\mathbf S}_2^*\cr
\vdots& \vdots& \ddots& \vdots\cr
{\mathbf S}_J^* & {\mathbf S}_J^* & \cdots& {\mathbf I}_G}
\right]
\left[\matrix{
{\mathbf m}_{1,j}^*\cr
{\mathbf m}_{2,j}^*\cr
\vdots\cr
{\mathbf m}_{J,j}^*} \right]
= \left[\matrix{
{\mathbf S}_1^*\cr
{\mathbf S}_2^*\cr
\vdots\cr
{\mathbf S}_J^*} \right] {\mathbf Y}_j^*,
\]
where ${\mathbf S}_j^*=({\mathbf I}_G-{\mathbf1}{\mathbf
1}^T/G){\mathbf S}_j$ is the centered smoother matrix, and $\mathbf
{S}_j$ is
defined the same as before.
The backfitting
estimators converge to the solution,
%
%
\begin{equation}
\left[\matrix{
\hat{{\mathbf m}}_{1,j}^*\cr
\hat{{\mathbf m}}_{2,j}^*\cr
\vdots\cr
\hat{{\mathbf m}}_{J,j}^*} \right]
= \left[\matrix{
{\mathbf I}_G & {\mathbf S}_1^* & \cdots& {\mathbf S}_1^*\cr
{\mathbf S}_2^* & {\mathbf I}_G & \cdots& {\mathbf S}_2^*\cr
\vdots& \vdots& \ddots& \vdots\cr
{\mathbf S}_J^* & {\mathbf S}_J^* & \cdots& {\mathbf I}_G}
\right]^{-1}
\left[\matrix{
{\mathbf S}_1^*\cr
{\mathbf S}_2^*\cr
\vdots\cr
{\mathbf S}_J^*} \right] {\mathbf Y}_j^*
\equiv{\mathbf M}^{-1}{\mathbf C}{\mathbf Y}_j^*,
\label{B2}
\end{equation}
provided that the inverse of ${\mathbf M}$ exists.

As in Opsomer (\citeyear{O00}), we define the additive smoother matrix as
\[
{\mathbf W}_k={\mathbf E}_k {\mathbf M}^{-1}{\mathbf C},
\]
where ${\mathbf E}_k$ is a partitioned matrix of dimension $G\times
GJ$ with an $G\times G$ identity matrix as the $k$th ``block'' and
zeros elsewhere. Thus the backfitting estimator for ${\mathbf
m}_{k,j}^*$ is
%
%
\begin{equation}\label{B3a}
\hat{\mathbf m}_{k,j}^*={\mathbf W}_k {\mathbf Y}_j^*.
\end{equation}
Denote by ${\mathbf m}_j^*=\sum_{k=1}^J{\mathbf m}_{k,j}^*$
and ${\mathbf W}_M=\sum_{k=1}^J{\mathbf W}_k$. The backfitting
estimator of ${\mathbf m}_j^*$ is then $\hat{\mathbf
m}_j^*={\mathbf W}_M
{\mathbf Y}_j^*$.
Let
${\mathbf W}_M^{[-k]}$ be the additive smoother matrix for the data
generated by the
$(J-1)$-variate regression model, $Y_{gj}'=\sum_{k'=1,\ne k}^J
m_{k',j}^*(X_{gk'})+\varepsilon_{gj}^*$.

If $\Vert{\mathbf S}_k^* {\mathbf W}_M^{[-k]}\Vert<1$ for
some $k\in
\{1,\ldots,J\}$, by Lemma 2.1 of
Opsomer (\citeyear{O00}), the backfitting estimators exist and are unique,
and
%
%
\begin{eqnarray} \label{B4}
{\mathbf W}_k&=&{\mathbf I}_G-\bigl({\mathbf I}_G-{\mathbf S}_k^*{\mathbf
W}_M^{[-k]}\bigr)^{-1}
({\mathbf I}_G-{\mathbf S}_k^*)\nonumber\\[-8pt]\\[-8pt]
&=&\bigl({\mathbf I}_G-{\mathbf S}_k^*{\mathbf W}_M^{[-k]}\bigr)^{-1}{\mathbf S}_k^*
\bigl({\mathbf I}_G-{\mathbf W}_M^{[-k]}\bigr).\nonumber
\end{eqnarray}
In this section we make the same assumption that is made in Opsomer
(\citeyear{O00}), that is, the inequality
$\Vert{\mathbf S}_k^* {\mathbf W}_M^{[-k]}\Vert<1$
holds.

For each $j$, $\hat{\mathbf m}_{k,j}^*$ estimates $\mathbf m_{k,j}^*$.
Define $\hat{\mathbf m}_{k,j}$ equals $-J\hat{\mathbf m}_{k,j}^*$ for
$k\ne j$ and
$J(J-1)^{-1}\hat{\mathbf m}_{k,j}^*$ for $k=j$.
Then $\hat{\mathbf m}_{k,j}$ estimates $\mathbf m_{k}$.
Since the variance of $\hat{\mathbf m}_{k,j}$ ($j\ne k$) is much bigger
than that of $\hat{\mathbf m}_{k,k}$,
taking the average over $j$ does not help reduce the variance of $\hat
{\mathbf m}_{k,k}$.
We will use $\hat{\mathbf m}_k\equiv\hat{\mathbf m}_{k,k}$ as an
estimate of $\mathbf m_{k}$.
The following theorem is a corollary of Theorem 3.1 in Opsomer
(\citeyear{O00}).
%
%
\begin{Theorem}\label{thbb}
Suppose that the conditions in Appendix \ref{appI} hold. If $X_{g1}$ is an
interior point, then as $G\to\infty$:
\begin{longlist}
\item The conditional bias of $\hat{m}_1(X_{g1})$ is
\begin{eqnarray*}
&&E\{\hat{m}_1(X_{g1})-m_1(X_{g1})|\mathbf{X}\}\\
&&\qquad=\mathbf e_g^T\bigl(I-\mathbf{S}_1^*\mathbf{W}_M^{[-1]}\bigr)^{-1}\\
&&\qquad\quad{}\times \biggl\{
\frac{\mu_2(K)}{2}h_1^2[{\mathcal D}^2\mathbf m_1-E(\mathbf
m_1'')]-\mathbf{S}
_1^*\mathbf{B}_{(-1)}
\biggr\}+o_p(h_1^2),
\end{eqnarray*}
where
$\mathbf{B}_{(-1)}=(\mathbf{W}_M^{[-1]}-\mathbf I_G)\mathbf m_{(-1)}$
and
$\mathbf m_{(-1)}=\sum_{k=2}^Jm_k$.
\item The conditional variance of $\hat{m}_1(X_{g1})$ is
\[
\operatorname{Var}\{\hat{m}_1(X_{g1})|\mathbf{X}\}=\frac{J}{J-1}\frac
{1}{Gh_1}\sigma
^2f_1^{-1}(X_{g1})\nu_0\bigl(K_{(1)}\bigr)+o_p\biggl(\frac{1}{Gh_1}\biggr).
\]
\end{longlist}
\end{Theorem}

As in Corollary 3.2 of Opsomer (\citeyear{O00}), if the covariates are
mutually independent, the conditional bias of $\hat{m}_1$ at an
interior observation point $X_{g1}$ is
\begin{eqnarray*}
&&E\{\hat{m}_1(X_{g1})-m_1(X_{g1})|\mathbf{X}\}
\\
&&\qquad= \frac{\mu_2(K_{(1)})}{2}h_1^2[m_1''(X_{g1})-E(m_1''(X_{g1}))]\\
&&\qquad\quad{}+O_p\bigl(1/\sqrt{G}\bigr)+o_p\Biggl(\sum_{j=1}^Jh_j^2\Biggr).
\end{eqnarray*}
This demonstrates that the estimators based on fitting
bivariate additive models and a multiple additive model
have the same asymptotic bias and variance in the interior points
when the covariates are independent. However, the estimator based on
fitting bivariate additive models is easy to implement.

\section{Simulations}\label{sec:d}

We here conduct simulations to compare the performance of the proposed
integration estimation method with the backfitting estimation.
To this end, we consider model (\ref{b1}) and
set $J=3$ and $G=3000$. The first variable $X_{g1}$ is generated from a
mixture distribution; that
is, $X_{g1}$ is simulated from the probability distribution $0.0004
\times(x-6)^3I(6 < x < 16)$ with
probability $0.6$ and from the uniform distribution over $[6, 16]$ with
probability $0.4$. The other two
variables $X_{gk}$ ($k=2,3$) are generated from model (\ref{b3}) with
$b_G=G^{-\gamma}$ and
$u_{gk}\sim_{\mathrm{i.i.d.}} N(0,1)$ where
$\gamma=0.05,0.1$ and $0.2$ are used to control the correlation
between $X_{gk}$ and $X_{g1}$. It is easy to calculate that $\gamma=
0.05, 0.1$ and 0.2 correspond to correlations 0.9919, 0.9962 and 0.9992
between $X_{g1}$ and $X_{g2}$, respectively. The correlations between
$X_{g2}$ and $X_{g3}$ are very close to the correlations between
$X_{g1}$ and $X_{g2}$ for different values of $\gamma$. The treatment
effect $\alpha_g$ is generated from the double exponential
distribution $\frac{1}{2}\exp(-|x|)$.
The response variable $Y_{kg}$ is simulated from model (\ref{b1}) with
$m_1(x)=\sqrt{5}(\sin(x)-0.2854)$, $m_2(x)=0.01(x-11)^3-0.2913$,
$m_3(x)=0.2\exp(x/5)-3.0648$ and
$\varepsilon_{kg}\sim_{\mathrm{i.i.d.}} N(0,1)$.

The mean square error (MSE) is employed to evaluate the performance of
different estimation methods.
The MSE of an estimate $\hat{m}_j$ of the function\vspace*{1pt} $m_j$ and the MSE of
an estimate $\hat{\bolds\alpha} = (\hat{\alpha}_1, \ldots, \hat{\alpha
}_G)^T$ of the vector $\bolds\alpha= (\alpha_1, \ldots, \alpha_G)^T$ are
defined, respectively, as follows:
\begin{eqnarray*}
\operatorname{MSE}(\hat{m}_j) &=& \frac{1}{G}\sum_{g=1}^G \bigl(\hat
{m}_j(X_{gj}) -
m_j(X_{gj}) \bigr)^2,\\
\operatorname{MSE}(\hat{\bolds\alpha}) &=& \frac{1}{G}\sum_{g=1}^G (\hat
{\alpha}_g -
\alpha_g)^2.
\end{eqnarray*}

%
\begin{figure}

\includegraphics{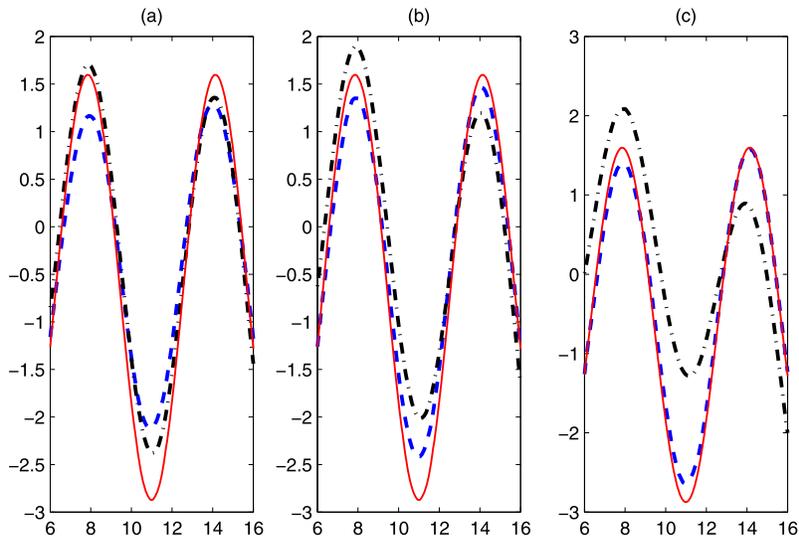}

\caption{Estimates of function
$m_1=\sqrt{5}(\sin(x)-0.2854)$. Solid curves (red): the true function.
Dashed curves (blue): the integration estimation method. Dash-dotted
curves (black): the pooled backfitting method. \textup{(a)} $\gamma=0.05$;
\textup{(b)} $\gamma=0.1$; \textup{(c)} $\gamma=0.2$.}\label{fig1}
\end{figure}

%
\begin{figure}[b]

\includegraphics{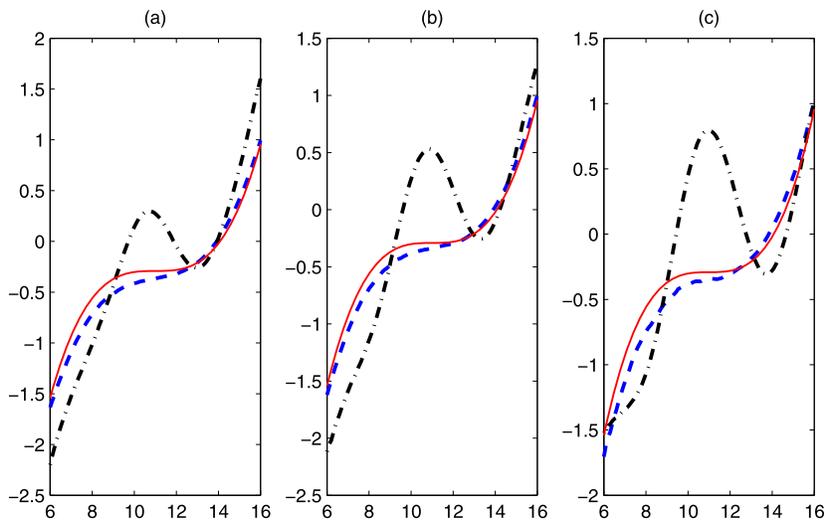}

\caption{The same as Figure \protect\ref{fig1} except that the
estimated function is
$m_2(x)=0.01(x-11)^3-0.2913$. \textup{(a)} $\gamma=0.05$; \textup{(b)} $\gamma=0.1$;
\textup{(c)} $\gamma=0.2$.}\label{fig2}
\end{figure}

The integration estimation procedure
and the pooled backfitting method are applied to estimate $m_1(\cdot)$
at 100 equispaced grid
points over the interval $[6,16]$ using $500$ simulated datasets.
For the backfitting method, we first tried the Gaussian kernel and the
optimal data-driven bandwidth rule in Opsomer and Ruppert (\citeyear
{OR98}) and
noticed that the estimated curves for the backfitting estimators were
over-smoothed when $\gamma$ is smaller. Following the reviewers'
suggestions, we then used a smaller bandwidth, that is, 0.4 times the
optimal bandwidth. For the integration method, its performance is not
sensitive to the choice of bandwidth, as long as it is not chosen too
large (see Theorem \ref{Th2}). Thus we just chose a reasonably small
one. The
medians of the fitted curves over 500 simulations are summarized in
Figure \ref{fig1}.
It is seen from Figure \ref{fig1} that,
when $\gamma$ becomes larger,
the correlation between covariates gets higher and the
backfitting method performs worse while the integration method becomes
better. In fact, when $\gamma=0.2$,
our integration procedure gives almost perfect estimates of the true
function: very little bias is involved. Similarly, we
estimate the functions $m_2(x)$ and $m_3(x)$. The estimated curves
are depicted in Figures \ref{fig2} and \ref{fig3}. It can be seen that
due to the high correlation, the pooled backfitting method gives
estimates that are highly biased while our integration method produces
almost perfect fits. The variations of the estimates are accessed by
MSE, and the median of these 500 MSEs can be found in Table \ref{tab1}.

%
\begin{figure}

\includegraphics{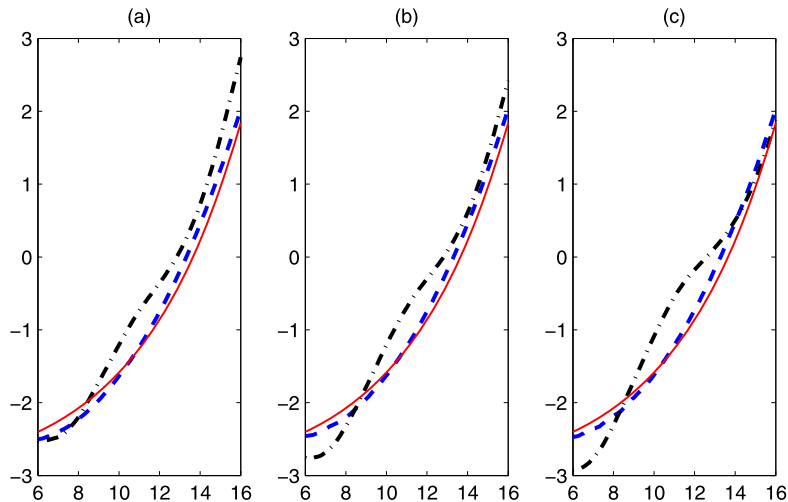}

\caption{The same as Figure \protect\ref{fig1} except that the
estimated function is $m_3(x)=0.2\exp(x/5)-3.0648$. \textup{(a)} $\gamma=0.05$;
\textup{(b)} $\gamma=0.1$; \textup{(c)} $\gamma=0.2$.}\label{fig3}
\end{figure}


%
\begin{table}[b]
\caption{Medians of MSEs for the estimated $m_j$ and $\bolds
\alpha$}\label{tab1}
\begin{tabular*}{\tablewidth}{@{\extracolsep{\fill}}lcccccc@{}}
\hline
& \multicolumn{3}{c}{\textbf{Integration estimation method}} &
\multicolumn{3}{c@{}}{\textbf{Pooled backfitting method}} \\[-4pt]
& \multicolumn{3}{c}{\hrulefill} &
\multicolumn{3}{c@{}}{\hrulefill} \\
& $\bolds{\gamma=0.05}$ & $\bolds{\gamma=0.1}$ & $\bolds{\gamma=0.2}$
& $\bolds{\gamma=0.05}$ & $\bolds{\gamma=0.1}$ & $\bolds{\gamma=0.2}$ \\
\hline
$m_1$ & 0.1471 & 0.0698 & 0.1032& 0.0774 & 0.2411& 0.8169 \\
$m_2$ & 0.0121 & 0.0177 & 0.0689 & 0.2310 & 0.1746 & 0.2343 \\
$m_3$ & 0.0202 & 0.0245 & 0.0750 & 0.1007 & 0.0754& 0.1254 \\
[4pt]
$\bolds\alpha$ & 0.3542 & 0.3565 & 0.3647 & 0.3963& 0.4125& 0.4962\\
\hline
\end{tabular*}
\end{table}

Now we estimate $\alpha_g, g=1,\ldots, G$. For each of the 500
simulated data sets, let $\hat{\alpha}_{gj} = Y_{gj} - \hat{m}_{j}(X_{gj})$,
for $g=1,\ldots, G$ and $j=1,\ldots, J$. Then for each of the simulated
data sets we estimate $\alpha_g$ as
\[
\hat{\alpha}_g = \frac{1}{J}\sum_{j=1}^J\hat{\alpha}_{gj}.
\]
The performance of $\hat{\bolds\alpha}$ is evaluated by MSE. The median of
the 500 MSEs is then calculated. Table \ref{tab1} reports the medians
of MSEs obtained by using the integration and pooled backfitting methods.
The integration estimation method dominates the backfitting method in
almost all cases.

\section{Real data example}\label{sec:e}

\subsection{Microarray data analysis}
We apply our new estimation methods to the Neuroblastoma data set
collected and analyzed by Fan et al. (\citeyear{Fan05a}). Neuroblastoma
is the
most frequent solid extra cranial neoplasia in children. Various
studies have suggested that microphage migration inhibitory factor
(MIF) may play an important role in the development of
neuroblastoma. To understand the impact of MIF reduction on
neuroblastoma cells, the global gene expression of the neuroblastoma
cell with MIF-suppressed is compared to those without MIF
suppression using Affymetrix GeneChips. Among extracted detection
signals, only genes with all detection signals greater than 50 were
considered, resulting in 13,980 genes in three control and treatment
arrays, respectively. The details of the design and experiments were
given by Fan et al. (\citeyear{Fan05a}).

For this DNA microarray data set, $J=3$ and $G=13\mbox{,}980$. Model
(\ref{b1}) was used in Fan et al. (\citeyear{Fan05a}) to assess the intensity
and treatment effects on genes with $m_j(\cdot)$ representing the
intensity effect for the $j$th array and $\alpha_g$ denoting the
treatment effect on gene $g$. As discussed in
Section \ref{sec:intro}, model (\ref{b1}) leads to the additive
model
%
%
\begin{equation}\label{eq001}
Y_g^{(k)}=m_1(X_{g1})-m_k(X_{gk})+\varepsilon_{g}^{(k)},\qquad
k=2,3,
\end{equation}
where $Y_{g}^{(k)}=Y_{g1}-Y_{gk}$. Now we
fit the data using model (\ref{eq001}) and estimate the components by
the integration and pooled backfitting methods. The resulting
estimates indicate similar forms of the intensity effects
for different slides, as presented in Figure \ref{fig4}. However, the
integration and pooled backfitting estimates differ substantially which
raises a question about which estimate is more reliable.


%
\begin{figure}[b]

\includegraphics{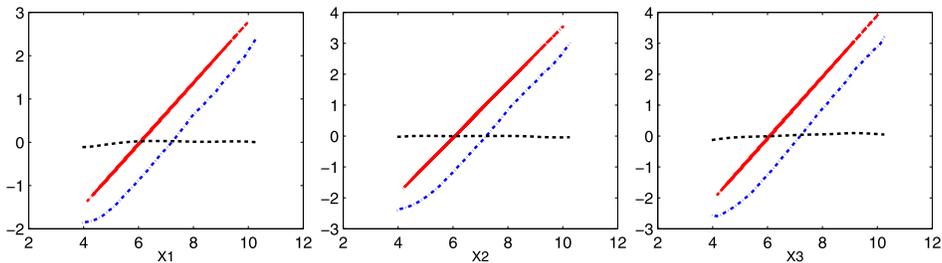}

\caption{Fitted regression curves as estimates of the
intensity effects for different arrays.
Left panel: for the first array, middle panel: for the second array,
right panel: for the third array;
dashed (black): the backfitting estimate, dashed-dotted (blue): the
integration estimate, solid (red): the linear regression.}\label{fig4}
\end{figure}

In the implementation of the backfitting method, we encounter an
almost singular matrix problem when using the Matlab software due to
the highly correlated log intensities $X_{gj}$, which leads to the
extremely low rate of convergence and unreliable results, even
though it reports the final estimates. Hence, by intuition and the
previous theory the integration estimation is better. In addition,
since both estimation methods lead to roughly linear forms of the
intensity effects functions $m_j(\cdot)$ for $j=1,2,3$, the linear
model seems
plausible. This suggests that we should fit the data using the linear model
\[
Y_{g}^{(k)}=\beta_0+\beta_1X_{g1}+\beta_2X_{gk}+\varepsilon_{g}^{(k)},
\]
as an alternative of model (\ref{eq001}).
For the integration estimation method
we do not have the singularity nor convergence problem. Thus we
still work with model (\ref{eq001}).


%
\begin{table}
\caption{The standard deviation of residuals from different
estimation for the additive models with~different~covariates}
\label{tab2}
\begin{tabular*}{\tablewidth}{@{\extracolsep{\fill}}lccc@{}}
\hline
\textbf{Covariates in models} & \textbf{Integration estimation} &
\textbf{Pooled backfitting} & \textbf{LSE} \\
\hline
$(X_1,X_2)$ & 0.430 & 0.451 & 0.429\\
$(X_1,X_3)$ & 0.421 & 0.462 & 0.428\\
$(X_2,X_3)$ & 0.375 & 0.455 & 0.455\\
\hline
\end{tabular*}
\end{table}

Figure \ref{fig4} displays
the estimated functions for each array. The pooled
backfitting estimates are almost flat and deviate far away from the
trends revealed by the least squares estimates (LSEs) for the linear model,
but the integration estimates share similar trends as the LSEs. It
seems that the estimated intensities from the integration method are
increasing and have a similar trend which suggests that the
intensity effects for the three slides are similar. Table \ref{tab2}
reports the standard deviation of residuals from the different
estimation methods. It favors the integration method.

\subsection{Interest rate data analysis}
In this subsection, we analyze the interest rates data introduced in
the \hyperref[sec:intro]{Introduction}. For simplicity, we consider the model with two
additive functions
\[
X_t = \mu+ m_1(X_{t-1}) + m_2 (X_{t-2}) + \varepsilon_t.
\]
%
Note that the above model is exactly model (\ref{b2}). Thus backfitting
and integration methods can be used to estimate the additive
components. Our integration method can be easily extended to the case
where there are three or more additive functions. Figure \ref{fig6}
shows the estimated functions $m_1(x)$ and $m_2(x)$ by using the
integration and pooled backfitting methods. Figure \ref{fig7} shows the
%
%
\begin{figure}

\includegraphics{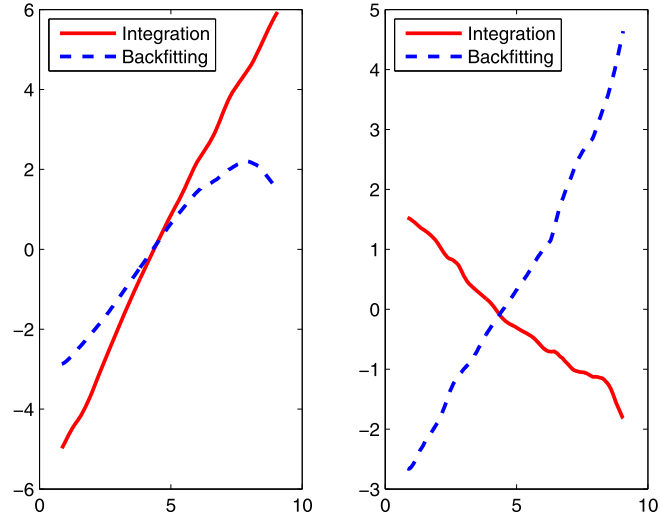}

\caption{Left panel: the estimated curve for $m_1(x)$;
right panel: the estimated curve for $m_2(x)$.}\label{fig6}
\end{figure}
%
%
\begin{figure}[b]

\includegraphics{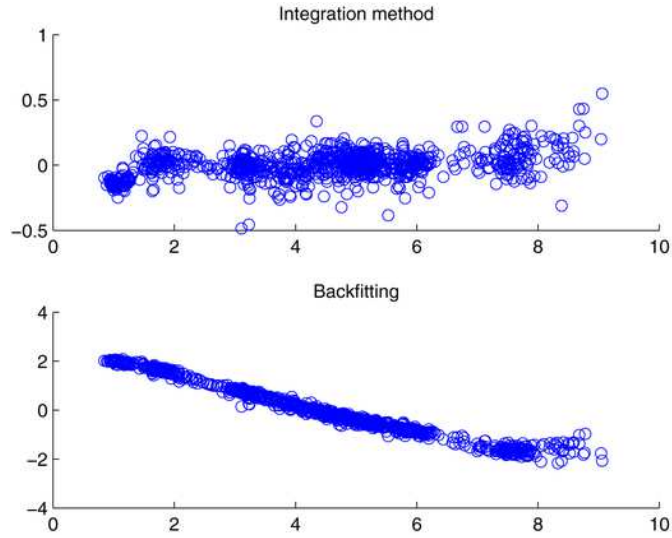}

\caption{Top panel: residual plot when the integration
method is used; bottom panel: residual plot when the backfitting method
is used.}\label{fig7}
\end{figure}
corresponding residuals,
which demonstrates that the integration method provides much better
fitting than the pooled backfitting method. Failure of the latter method
is the result of highly correlated covariates [see also Figure \ref
{fig0}(right)] in the fitted model.

\section{Discussion}

In this article we have proposed several estimation methods for
additive models
when its covariates are highly correlated and nonhighly correlated.
We derived asymptotic normality of the proposed estimators
and illustrated their performance in finite samples via simulations.
The performance of the proposed methodology was also demonstrated by
two real data examples.

Many problems remain open for the array-dependent model.
Examples include:
\begin{longlist}
\item Investigation of the asymptotic normality of the backfitting
estimators when the covariates are highly correlated.

\item Establishing the asymptotic distribution of the estimators
in (\ref{eq1a}).

\item Test if the nonparametric functions $m_j$ have certain
parametric forms.
The generalized likelihood ratio tests can be used [see Fan and
Jiang (\citeyear{FJ05}, \citeyear{FJ07})].
\end{longlist}

\begin{appendix}
\section{Conditions}\label{appI}

\begin{longlist}
\item The kernel $K(\cdot)$ is a continuous and symmetric function
and has compact support,
and its first derivatives had a finite number of sign changes over its support.

\item The densities of $f_j$'s are bounded and continuous, have
compact support and their first derivatives
have a finite number of sign changes over their supports.
Also, $f_j(x_j)>0$ for all $x_j\in\operatorname{supp}(f_j)$.

\item As $G\to\infty$, $h_j\to0$, $h\to0$,
$Gh_j/\log G\to\infty$ and $Gh/\log G\to\infty$.

\item The second derivatives of $m_j$ exist and are continuous
and bounded.
\end{longlist}

\section{Proofs of theorems}\label{appII}

\begin{pf*}{Proof of Theorem \protect\ref{Th1}}
Let $\mathbf{Y}_{(k)}=(Y_1^{(k)},\ldots,Y_G^{(k)})^T$,
\[
\tilde{\bolds\varepsilon}^{(k)}=\bigl(\tilde{\varepsilon}^{(k)}_1,\ldots
,\tilde
{\varepsilon}^{(k)}_G\bigr)^T,
\]
$\tilde{m}_g(X_{gk})=m_{k1}(X_{gk})+b_Gu_{gk}m_1'(X_{gk})$
and
$\tilde{\mathbf m}=(\tilde{m}_1(X_{1k}),\ldots,\tilde{m}_G(X_{Gk}))^T$.
Then (\ref{b2a}) becomes
\[
\mathbf{Y}_{(k)}=\tilde{\mathbf m}+\tilde{\bolds\varepsilon}^{(k)}.
\]
By (\ref{b7aa}), we have
%
%
\begin{eqnarray}\label{ap1}\quad
\hat{\theta}(x)-\theta(x)&=&(\mathbf{Z}^T\mathbf K\mathbf{Z})^{-1}\mathbf{Z}^{T}
\mathbf K\tilde{\bolds\varepsilon}^{(k)}+(\mathbf{Z}^T\mathbf K\mathbf
{Z})^{-1}\mathbf{Z}^{T}\mathbf K[\tilde
{\mathbf m}-\mathbf{Z}\theta(x)]\nonumber\\[-8pt]\\[-8pt]
&=&\mathbf{V}(x)+\mathbf{B}(x).\nonumber
\end{eqnarray}
Note that for $|X_{gk}-x|\leq h$,
\begin{eqnarray*}
\tilde{m}_g(X_{gk})
&=&m_{k1}(x)+m'_{k1}(x)(X_{gk}-x)\\
&&{}+\tfrac{1}{2}m''_{k1}(x)(X_{gk}-x)^2+o(X_{gk}-x)^2 \\
&&{}+b_Gu_{gk}\bigl\{m_1'(x)+m_1''(x)(X_{gk}-x)\\
&&\hspace*{55.12pt}{}+\tfrac{1}{2}m_1^{(3)}(x)(X_{gk}-x)^2\bigr\}+o(X_{gk}-x)^2b_G\\
&=&Z_g^T\theta(x)+\tfrac{1}{2}m''_{k1}(x)(X_{gk}-x)^2\\
&&{}+\tfrac{1}{2}m_1^{(3)}(x)b_Gu_{gk}(X_{gk}-x)^2+o(h^2+hb_G),
\end{eqnarray*}
uniformly for $g=1,\ldots,G$.
Let
\[
\mathbf{X}= \left[\matrix{
(X_{1k}-x)^2h^{-2}& (X_{1k}-x)^2h^{-2}u_{1k}\cr
\vdots&\vdots\cr
(X_{Gk}-x)^2h^{-2}& (X_{Gk}-x)^2h^{-2}u_{Gk}
} \right].
\]
Then
%
%
\begin{equation}\label{newmy}
\tilde{\mathbf m}_{1k}-\mathbf{Z}\theta(x)=\frac{h^2}{2}\mathbf
{X}\left[\matrix{
m''_{k1}(x)\cr
m_1^{(3)}(x)b_G}
\right]+o({\mathbf1})(h^2+b_Gh^2),
\end{equation}
where
${\mathbf1}$ is a $G\times1$ vector with all elements being $1$'s,
and hence
%
%
\begin{equation}\label{ap2}\quad
\mathbf{B}(x)=(\mathbf{Z}^T\mathbf K\mathbf{Z})^{-1}\mathbf{Z}^T\mathbf K\left\{\mathbf
{X}\left[\matrix{
\dfrac{h^2}{2}m''_{k1}(x)\vspace*{2pt}\cr
hb_Gm''_1(x)}
\right]+o({\mathbf1})(h^2+h^2b_G) \right\}.
\end{equation}
Let $\mathbf{S}_T=\mathbf{Z}^T\mathbf K\mathbf{Z}$.
Then
$\mathbf{S}_T=\sum_{g=1}^GK_h(X_{gk}-x)Z_gZ_g'$, and
the $(i,j)$th element of $\mathbf{S}_T$
is
\[
\mathbf{S}_{T,ij}=\cases{
\displaystyle\sum_{g=1}^GK_h(X_{gk}-x)(X_{gk}-x)^{i+j-2},&\quad for $1\leq i,j\leq2$;\cr
\displaystyle\sum_{g=1}^GK_h(X_{gk}-x)(X_{gk}-x)^{i+j-4}b_Gu_{gk},&\quad for $i=1,2;
j=3,4$;\cr
\displaystyle\sum_{g=1}^GK_h(X_{gk}-x)(X_{gk}-x)^{i+j-4}b_Gu_{gk},&\quad for
$i=3,4;j=1,2$;\cr
\displaystyle\sum_{g=1}^GK_h(X_{gk}-x)(X_{gk}-x)^{i+j-6}b_G^2u_{gk}^2, &\quad for
$i,j=3,4$.}
\]
Directly computing the mean and variance, we obtain that:
\begin{longlist}
\item for $1\leq
i,j\leq2$,
\begin{eqnarray*}
G^{-1}h^{-(i+j-2)}\mathbf{S}_{T,ij}&=&G^{-1}\sum
_{g=1}^GK_h(X_{gk}-x)(X_{gk}-x)^{i+j-2}h^{-(i+j-2)}\\
&=&f_k(x)\mu_{i+j-2}(K)+O_p\bigl(h+1/\sqrt{Gh}\bigr) ;
\end{eqnarray*}
\item for $i=1,2$ and
$j=3,4$, or $i=3,4$ and $j=1,2$,
\begin{eqnarray*}
&&G^{-1}h^{-(i+j-4)}b_G^{-1}\mathbf{S}_{T,ij}\\
&&\qquad=G^{-1}\sum_{g=1}^G
K_h(X_{gk}-x)(X_{gk}-x)^{i+j-4}h^{-(i+j-4)}u_{gk}\\
&&\qquad=O_p\bigl(1/\sqrt{Gh}\bigr);
\end{eqnarray*}
\item for $i,j=3,4$,
\begin{eqnarray*}
&&G^{-1}h^{-(i+j-6)}b_G^{-2}\mathbf{S}_{T,ij}
\\
&&\qquad=G^{-1}\sum_{g=1}^GK_h(X_{gk}-x)(X_{gk}-x)^{i+j-6}h^{-(i+j-6)}
u_{gk}^2\\
&&\qquad=f_k(x)\mu_{i+j-6}(K)+O_p\bigl(h+1/\sqrt{Gh}\bigr).
\end{eqnarray*}
\end{longlist}
Therefore,
%
%
\begin{equation}\label{newmy1}
G^{-1}\mathbf{H}^{-1}\mathbf{S}_T\mathbf{H}^{-1}=f_k(x)\mathbf
{S}+O_p({\mathbf1}{\mathbf
1}^T)\biggl(h+\frac{1}{\sqrt{Gh}}\biggr).
\end{equation}
By simple algebra and (\ref{newmy}), we have
\begin{eqnarray*}
&&G^{-1}\mathbf{H}^{-1}\mathbf{Z}^T\mathbf K[\tilde{\mathbf m}_{1k}-\mathbf
{Z}\theta(x)]
\\
&&\qquad= G^{-1}\mathbf{H}^{-1}\mathbf{Z}^T\mathbf K\\
&&\qquad\quad{}\times\left\{\mathbf{X}\left[\matrix{
\dfrac{h^2}{2}m''_{k1}(x)\vspace*{2pt}\cr
\dfrac{h^2}{2}m_1^{(3)}(x)b_G
} \right]+o({\mathbf1})(h^2+b_Gh) \right\}\\
&&\qquad=\mathbf{A}\left[\matrix{
\dfrac{h^2}{2}m''_{k1}(x)\vspace*{2pt}\cr
hb_Gm''_1(x)}
\right]+o(h^2+b_Gh),
\end{eqnarray*}
uniformly for components
where $\mathbf{A}=(A_{ij})$ is a $4\times2$ matrix with
\[
A_{ij}=\cases{
\displaystyle G^{-1}\sum_{g=1}^G K_h(X_{gk}-x)h^{-(i+j)}(X_{gk}-x)^{i+j}, &\quad for
$i=1,2$ and
$j=1$;\cr
\displaystyle G^{-1}\sum_{g=1}^G K_h(X_{gk}-x)h^{-i}(X_{gk}-x)^{i}u_{gk}, &\quad for
$i=1,2$ and
$j=2$;\cr
\displaystyle G^{-1}\sum_{g=1}^G K_h(X_{gk}-x)h^{-2}(X_{gk}-x)^2u_{gk}, &\quad for
$i=3,j=1$;\cr
\displaystyle G^{-1}\sum_{g=1}^G K_h(X_{gk}-x)h^{-1}(X_{gk}-x)u_{gk}^2,&\quad for
$i=3,j=2$;\cr
\displaystyle G^{-1}\sum_{g=1}^G K_h(X_{gk}-x)h^{-3}(X_{gk}-x)^3u_{gk},&\quad for
$i=4,j=1$;\cr
\displaystyle G^{-1}\sum_{g=1}^G K_h(X_{gk}-x)h^{-3}(X_{gk}-x)^3u_{gk}^2,&\quad for
$i=4,j=2$.}
\]
Directly computing the mean and variance of $A_{ij}$, we obtain that
$\mathbf{A}=f_k(x)\mathbf{C}+o(h^2+b_Gh)$, uniformly for components.
Then
\begin{eqnarray*}
&&G^{-1}\mathbf{H}^{-1}\mathbf{Z}^T\mathbf K[\tilde{\mathbf m}_{1k}-\mathbf
{Z}\theta(x)]\\
&&\qquad=f_k(x) \left[\matrix{
\mu_2&0\cr
\mu_3&0\cr
0&\mu_2\cr
0&\mu_3}
\right] \left[\matrix{
\dfrac{h^2}{2}m''_{k1}(x)\vspace*{2pt}\cr
\dfrac{h^2}{2}b_Gm_1^{(3)}(x)}
\right]+o(h^2+b_Gh)\\
&&\qquad=
\frac{h^2}{2}f_k(x)\mathbf{C}
\bigl(m''_{k1}(x),m_1^{(3)}(x)b_G\bigr)^T+o(h^2+b_Gh),
\end{eqnarray*}
uniformly for components.
Thus
\begin{eqnarray*}
\mathbf{H}\mathbf{B}(x)&=&\mathbf{H}\mathbf{S}_T^{-1}\mathbf{Z}^T\mathbf K
[\tilde{\mathbf m}_{1k}-\mathbf{Z}\theta(x)]\\
&=&
(\mathbf{H}^{-1}\mathbf{S}_T\mathbf{H}^{-1})^{-1}\mathbf{H}^{-1}\mathbf
{Z}^T\mathbf K[\tilde{\mathbf m}_{1k}-\mathbf{Z}
\theta(x)]\\
&=&\frac{h^2}{2}\mathbf{S}^{-1}\mathbf{C}
\bigl(m''_{k1}(x),m_1^{(3)}(x)b_G\bigr)^T\bigl(1+o_p(1)\bigr).
\end{eqnarray*}
This combined with (\ref{ap1}) yields that
%
%
\begin{equation} \label{ap3}\hspace*{35pt}
\mathbf{H}\bigl(\hat{\theta}(x)-\theta(x)\bigr)-\frac{h^2}{2}\mathbf
{S}^{-1}\mathbf{C}
\bigl(m''_{k1}(x),b_Gm_1^{(3)}(x)\bigr)^T\bigl(1+o_p(1)\bigr)
=\mathbf{H}\mathbf{V}(x).
\end{equation}
By the definition of $\mathbf{V}(x)$ we have
\[
\mathbf{H}
\mathbf{V}(x)=\mathbf{H}\mathbf{S}_T^{-1}\mathbf{Z}^T\mathbf K\tilde{\bolds
\varepsilon}^{(k)}
=(\mathbf{H}^{-1}\mathbf{S}_T\mathbf{H}^{-1})^{-1}\mathbf{H}^{-1}\mathbf
{Z}^T\mathbf K\tilde{\bolds\varepsilon}{}^{(k)}.
\]
Plugging (\ref{newmy1}) into the right-hand side above, we establish that
\begin{eqnarray*}
\mathbf{H}\mathbf{V}(x)
&=&G^{-1}(f_k(x)\mathbf{S})^{-1}
\left[\matrix{
1&\cdots&1\vspace*{2pt}\cr
\dfrac{X_{1k}-x}{h}&\cdots&\dfrac{X_{Gk}-x}{h}\vspace*{3pt}\cr
u_{1k}&\cdots& u_{gk}\vspace*{0pt}\cr
\dfrac{X_{1k}-x}{h}u_{1k}&\cdots&\dfrac{X_{Gk}-x}{h}u_{Gk}}
\right]\mathbf K\left[\matrix{
\tilde{\varepsilon}_1^{(k)}\cr
\vdots\cr
\tilde{\varepsilon}_G^{(k)}}
\right]\\
&=&f_k^{-1}(x)\mathbf{S}^{-1}J_G(x),
\end{eqnarray*}
where
\[
J_G(x)= \left[\matrix{
\displaystyle G^{-1}\sum_{g=1}^GK_h(X_{gk}-x)\tilde{\varepsilon}_g^{(k)}\cr
\displaystyle G^{-1}\sum_{g=1}^GK_h(X_{gk}-x)\biggl(\frac{X_{gk}-x}{h}\biggr)\tilde{\varepsilon
}_g^{(k)}\cr
\displaystyle G^{-1}\sum_{g=1}^GK_h(X_{gk}-x)u_{gk}\tilde{\varepsilon}_g^{(k)}\cr
\displaystyle G^{-1}\sum_{g=1}^GK_h(X_{gk}-x)\biggl(\frac{X_{gk}-x}{h}\biggr)u_{gk}\tilde
{\varepsilon}_g^{(k)}}
\right].
\]
Under the working model (\ref{b3}), we obtain from (\ref{b2a}) that
%
%
\begin{equation}\label{newc}
\tilde{\varepsilon}_g^{(k)}=\tfrac
{1}{2}m_1''(X_{gk})b_G^2u_{gk}^2\bigl(1+o_p(1)\bigr)+\varepsilon_g^{(k)}.
\end{equation}
Using an argument similar to that for Lemma 7.3 of Jiang and Mack
(\citeyear{JM01}), we can show that
\begin{eqnarray*}
&&\sqrt{Gh}\bigl[\mathbf{H}\mathbf{V}(x)-\tfrac{1}{2}b_G^2m_1''(x)\mathbf
{S}^{-1}\mathbf{c}
^*\bigl(1+o_p(1)\bigr)+O_p\bigl(b_G^2/\sqrt{Gh}\bigr)\bigr]\\
&&\qquad\stackrel{\mathcal D}{\longrightarrow}
N(0,2f_k^{-1}(x)\sigma^2\mathbf{S}^{-1}\mathbf{V}\mathbf{S}^{-1}),
\end{eqnarray*}
which together with (\ref{ap3}) and
$f_1(x)=f_k(x)(1+o(1))$
leads to the result of the theorem.
\end{pf*}
\begin{pf*}{Proof of Corollary \protect\ref{col1}}
Let $e_3=(0,0,1,0)^T$.
Then
\begin{eqnarray*}
&&\hat{m}_1'(x;k)-m_1'(x)\\
&&\qquad=e_3^T\bigl(\hat{\theta}(x)-\theta(x)\bigr)\\
&&\qquad=e_3^T\frac{h^2}{2}\mathbf{H}^{-1}\mathbf{S}^{-1}\mathbf{C}
\bigl(m''_{k1}(x),b_Gm_1^{(3)}(x)\bigr)^T\bigl(1+o_p(1)\bigr)\\
&&\qquad\quad{} +e_3^Tf_k^{-1}(x)\mathbf{H}^{-1}\mathbf{S}^{-1}J_G(x).
\end{eqnarray*}
It is easy to verify that
$e_3^T\mathbf{H}^{-1}\mathbf{S}^{-1}=(0,0,b_G^{-1}\mu_0^{-1}(K),0)$. Then
%
%
\begin{eqnarray}\label{eqad}
&&\hat{m}_1'(x;k)-m_1'(x)\nonumber\\
&&\qquad=\frac{h^2}{2}\mu_2(K)\mu_0^{-1}(K)m_1^{(3)}(x)\bigl(1+o_p(1)\bigr)\\
&&\qquad\quad{}+f_k^{-1}(x)\mu_0^{-1}(K)b_G^{-1}G^{-1}\sum
_{g=1}^GK_h(X_{gk}-x)u_{gk}\tilde{\varepsilon}_g^{(k)}.\nonumber
\end{eqnarray}
This combined with the asymptotic normality of $J_G(x)$ completes the
proof of the corollary.
\end{pf*}
\begin{pf*}{Proof of Theorem \protect\ref{Th1a}}
By (\ref{eqad}),
\begin{eqnarray*}
&&\hat{m}_1'(x)-m_1'(x)\\
&&\qquad=\frac{h^2}{2}\mu_2(K)\mu_0^{-1}(K)m_1^{(3)}(x)\bigl(1+o_p(1)\bigr)\\
&&\qquad\quad{}+\frac{1}{J-1}\sum_{k=2}^J
f_k^{-1}(x)\mu_0^{-1}(K)b_G^{-1}G^{-1}\sum
_{g=1}^GK_h(X_{gk}-x)u_{gk}\tilde{\varepsilon}_g^{(k)}.%
\end{eqnarray*}
Then using (\ref{newc}), we obtain that
%
%
\begin{eqnarray} \label{eqad1}\qquad
&&\hat{m}_1'(x)-m_1'(x)\nonumber\\
&&\qquad=\frac{h^2}{2}\mu_2(K)\mu_0^{-1}(K)m_1^{(3)}(x)\bigl(1+o_p(1)\bigr)
\nonumber\\[-8pt]\\[-8pt]
&&\qquad\quad{} +\frac{1}{2}b_G m_1''(x)E(u_{1k}^3)\bigl(1+o_p(1)\bigr)\nonumber\\
&&\qquad\quad{}+\frac{1}{J-1}\sum_{k=2}^J
f_k^{-1}(x)\mu_0^{-1}(K)b_G^{-1}G^{-1}\sum
_{g=1}^GK_h(X_{gk}-x)u_{gk}{\varepsilon}_g^{(k)}.\nonumber
\end{eqnarray}
Let $B_G(x)=\frac{1}{J-1}\sum_{k=2}^J
f_k^{-1}(x)\mu_0^{-1}(K)b_G^{-1}G^{-1}\sum
_{g=1}^GK_h(X_{gk}-x)u_{gk}{\varepsilon}_g^{(k)}$. Then $E[B_G(x)]=0$.
Note that $E(u_{gk_1}u_{gk_2})=\rho(k_1,k_2)$ and
$E\{\varepsilon_g^{(k_1)}\varepsilon_g^{(k_2)}\}=\sigma^2$ for $k_1\ne
k_2$ and $2\sigma^2$ for $k_1= k_2$. It follows that
\begin{eqnarray*}
&&E\bigl[\sqrt{Gh}b_GB_G(x)\bigr]^2\\
&&\qquad=
\frac{1}{(J-1)^2}\sum_{k_1,k_2=2}^J f_{k_1}^{-1}(x)f_{k_2}^{-1}(x)
\mu_0^{-2}(K)\nonumber\\
&&\qquad\quad\hspace*{70.2pt}{}\times
E\bigl\{hK_h(X_{gk_1}-x)K_h(X_{gk_2}-x)\\
&&\hspace*{135.8pt}{}\times E\bigl(\varepsilon_g^{(k_1)}\varepsilon_g^{(k_2)}\bigr)E(u_{gk_1}u_{gk_2})\bigr\}
\nonumber\\
&&\qquad= \frac{2}{(J-1)^2}\sum_{k=2}^J f_{k}^{-2}(x)
\mu_0^{-2}(K)E\{hK_h^2(X_{gk}-x)\}\sigma^2\rho(k,k)\nonumber\\
&&\qquad\quad{} + \frac{1}{(J-1)^2}\sum_{k_1\ne k_2} f_{k_1}^{-1}(x)f_{k_2}^{-1}(x)
\mu_0^{-2}(K)\nonumber\\
&&\qquad\quad\hspace*{76.05pt}{}\times E\{hK_h(X_{gk_1}-x)K_h(X_{gk_2}-x)\}
\sigma^2\rho(k_1,k_2).
\end{eqnarray*}
Using $f_k(x)=f_1(x)(1+o(1))$ and
\begin{eqnarray*}
E\{hK_h(X_{gk_1}-x)K_h(X_{gk_2}-x)\}
&=&E\{hK_h^2(X_{g1}-x)\}\bigl(1+o(1)\bigr)\\
&=&f_1(x)\nu_0(K)\bigl(1+o(1)\bigr),
\end{eqnarray*}
we arrive at
\[
E\bigl[\sqrt{Gh}b_GB_G(x)\bigr]^2
=\rho\sigma^2 f_1^{-1}(x)\mu_0^{-2}(K)\nu_0(K)\bigl(1+o(1)\bigr),
\]
where
$\rho=\frac{1}{(J-1)^2}[\sum_{k=2}^J\rho(k,k)+\sum_{k_1=2}^J\sum
_{k_2=2}^J\rho(k_1,k_2)]$.
Therefore, $\sqrt{Gh}b_G\times\break B_G(x)$ is asymptotically normal with mean zero
and variance $\sigma_2^2(x)$.
This together with (\ref{eqad1}) completes the proof of the theorem.
\end{pf*}
\begin{pf*}{Proof of Theorem \protect\ref{Th2}}
Observing that
\begin{eqnarray*}
\hat{m}_1(x)
&=&-G^{-1}\sum_{g=1}^G\int_{x_0}^{X_{g1}}\hat{m}_1'(t) \,dt+\int
_{x_0}^x\hat{m}_1'(t) \,dt\\
&=&-G^{-1}\sum_{g=1}^G\int_{x_0}^{X_{g1}}[\hat{m}_1'(t)-m_1'(t)] \,dt+\int
_{x_0}^x[\hat{m}_1'(t)-m_1'(t)] \,dt\\
&&{} + m_1(x)-G^{-1}\sum_{g=1}^Gm_1(X_{g1})
\end{eqnarray*}
and
$G^{-1}\sum_{g=1}^Gm_1(X_{g1})=O_p(G^{-1/2})$,
we obtain that
\begin{eqnarray*}
&&\hat{m}_1(x)-m_1(x)\\
&&\qquad=-G^{-1}\sum_{g=1}^G\int_{x_0}^{X_{g1}}[\hat{m}_1'(t)-m_1'(t)] \,dt\\
&&\qquad\quad{} +\int_{x_0}^x[\hat{m}_1'(t)-m_1'(t)] \,dt+O_p(G^{-1/2})\\
&&\qquad=G^{-1}\sum_{g=1}^G\int_{X_{g1}}^x[\hat{m}_1'(t)-m_1'(t)] \,dt+O_p(G^{-1/2}).
\end{eqnarray*}
%
Let
$J_{k,G}=G^{-1}\sum_{g=1}^Gu_{gk}\varepsilon_g^{(k)}\int
_{X_{g1}}^xK_h(X_{gk}-t) \,dt$.
Then by (\ref{eqad1}) and simple algebra,
\begin{eqnarray*}
&&\hat{m}_1(x)-m_1(x)\\
&&\qquad=\frac{h^2}{2}\mu_2(K)\mu
_0^{-1}(K)[m_1''(x)-Em_1''(X_{11})]\bigl(1+o_p(1)\bigr)\nonumber\\
&&\qquad\quad{}+\frac{1}{2}b_GE(u_{1k}^3)[m_1'(x)-Em_1'(X_{11})]\bigl(1+o_p(1)\bigr)\nonumber\\
&&\qquad\quad{}+O_p(G^{-1/2})+\mu_0^{-1}(K)b_G^{-1}(J-1)^{-1}\sum_{k=2}^Jf_k^{-1}(x)J_{k,G}.
\end{eqnarray*}
Let
\begin{eqnarray*}
C_G(x)&=&\frac{1}{J-1}\sum_{k=2}^J
f_k^{-1}(x)\mu_0^{-1}(K)b_G^{-1}G^{-1}\\
&&{}\times\sum_{g=1}^Gu_{gk}{\varepsilon
}_g^{(k)}\int_{X_{g1}}^xK_h(X_{gk}-t)\, dt.
\end{eqnarray*}
Then
\begin{eqnarray*}
&&E\bigl[\sqrt{G}b_GC_G(x)\bigr]^2\\
&&\qquad= \frac{2}{(J-1)^2}\sum_{k=2}^J f_{k}^{-2}(x)
\mu_0^{-2}(K)E \biggl\{\int_{X_{g1}}^x K_h(X_{gk}-t) \,dt \biggr\}^2\sigma^2\rho
(k,k)\nonumber\\
&&\qquad\quad{}+\frac{1}{(J-1)^2}\sum_{k_1\ne k_2} f_{k_1}^{-1}(x)f_{k_2}^{-1}(x)
\mu_0^{-2}(K)\nonumber\\
&&\qquad\quad\hspace*{76.2pt}{}\times E \biggl\{\int_{X_{g1}}^xK_h(X_{gk_1}-t) \,dt\\
&&\qquad\quad\hspace*{76.2pt}\hspace*{24.6pt}{}\times \int
_{X_{g1}}^xK_h(X_{gk_2}-t) \,dt \biggr\}
\sigma^2\rho(k_1,k_2).
\end{eqnarray*}
Since 
\begin{eqnarray*}
&&\hspace*{30.75pt}E \biggl\{\int_{X_{g1}}^xK_h(X_{gk_1}-t) \,dt\int_{X_{g1}}^xK_h(X_{gk_2}-t)
\,dt \biggr\}\\
&&\hspace*{30.75pt}\qquad=E \biggl\{\int_{X_{g1}}^xK_h(X_{g1}-t) \,dt \biggr\}^2\bigl(1+o(1)\bigr)\\
&&\hspace*{30.75pt}\qquad= \int_{-\infty}^{\infty}f_1(u) \,du \biggl\{\int_u^xK_h(u-t) \,dt \biggr\}^2
\\
&&\hspace*{30.75pt}\qquad=\frac{1}{4}\mu_0^2(K)\bigl(1+o(1)\bigr),
\\
&&\hspace*{-30.75pt}E\bigl[\sqrt{G}b_GC_G(x)\bigr]^2
=\frac{1}{4}f_1^{-2}(x)\sigma^2\rho+o(1).
\end{eqnarray*}
Therefore, $\sqrt{G}b_GC_G(x)$ is asymptotically normal with mean zero
and variance $\sigma^2(x)$.
\end{pf*}
%
%
\begin{pf*}{Proof of Theorem \protect\ref{thb}}
As in Opsomer and Ruppert (\citeyear{OR97}), we let
\[
\mathbf Q_{m_1}(x_1)= \left[\matrix{
(X_{11}-x_1)^2\cr
\vdots\cr
(X_{G1}-x_1)^2}
\right]
\frac{\partial m_1(x_1)}{\partial x_1^2},\qquad
\mathbf Q_1= \left[\matrix{
\mathbf s^T_{1,X_{11}}\mathbf Q_{m_1}(\mathbf{X}_{11})\cr
\vdots\cr
\mathbf s^T_{1,X_{G1}}\mathbf Q_{m_1}(\mathbf{X}_{G1})}
\right],
\]
and similarly for $\mathbf Q_{m_k}(x_k)$ and $\mathbf Q_k$.
Let $h_j^2=h_j^2{\mathbf1}$. Then by the proof of Theorem~4.1 of
Opsomer and Ruppert (\citeyear{OR97}),
\[
(\mathbf I_G-\mathbf{S}_1^*\mathbf{S}_k^*)^{-1}(\mathbf I_G-\mathbf
{S}_1^*)\mathbf m_k
=\mathbf m_k+\tfrac{1}{2}(\mathbf I_G-\mathbf{S}_1^*\mathbf
{S}_k^*)^{-1}\mathbf{S}_1^*\mathbf
Q_k+o_p(\mathbf h_k^2)
\]
and
\[
(\mathbf I_G-\mathbf{S}_1^*\mathbf{S}_k^*)^{-1}(\mathbf I_G-\mathbf
{S}_1^*)\mathbf m_1
=\bar{\mathbf m}_1-\tfrac{1}{2}(\mathbf I_G-\mathbf{S}_1^*\mathbf
{S}_k^*)^{-1}\mathbf
Q_1^*+o_p(\mathbf h_1^2),
\]
where
$\mathbf Q_1^*=(\mathbf I_G-\mathbf1\mathbf1^T/G)\mathbf Q_1$.
Thus,
%
%
\begin{equation}\label{p1}
E(\hat{\mathbf m}_1-\mathbf m_1 |\mathbf{X})=\tfrac{1}{2}(\mathbf
I_G-\mathbf{S}_k^*\mathbf{S}
_1^*)^{-1}(\mathbf Q_1^*-\mathbf{S}_1^*\mathbf Q_k)+o_p(\mathbf
h_1^2+\mathbf h_k^2).
\end{equation}
By (\ref{B3}),
$\hat{\mathbf m}_1^{(j)}=\mathbf{W}_{1j}\mathbf m+\mathbf
{W}_{1j}\varepsilon_g^{(j)}$.
Note that
$\operatorname{Var}(\varepsilon_g^{(j)} |\mathbf{X})=2\sigma^2I_G$
and for $j\ne k$,
$\operatorname{Cov}(\varepsilon_g^{(j)},\varepsilon_g^{(k)} |\mathbf
{X})=\sigma^2\mathbf I_G$.
It follows that
%
%
\begin{eqnarray}\hspace*{32pt}
&&\operatorname{Cov}\bigl\{\hat{m}_1^{(j)}(X_{g1}),\hat{m}_1^{(k)}(X_{g1})
|\mathbf{X}\bigr\}\nonumber
\\
&&\qquad= \sigma^2\{1-\mathbf e_g^T(\mathbf I-\mathbf{S}_1^*\mathbf{S}_j^*)^{-1}
(\mathbf I-\mathbf{S}_1^*)\mathbf e_g
-\mathbf e_g^T(\mathbf I-\mathbf{S}_1^*\mathbf{S}_k^*)^{-1}(\mathbf
I-\mathbf{S}_1^*)\mathbf e_g\\
&&\qquad\quad\hspace*{46.4pt}{} +\mathbf e_g^T(\mathbf I-\mathbf{S}_1^*\mathbf{S}_j^*)^{-1}
(\mathbf I-\mathbf{S}_1^*)(\mathbf I-\mathbf{S}_1^*)^T(\mathbf I-\mathbf
{S}_1^*\mathbf{S}_k^*)^{-T}\mathbf e_g\}\nonumber
\end{eqnarray}
and
%
%
\begin{eqnarray} \label{eqnew1}\quad
&&\operatorname{Var}\bigl\{\hat{m}_1^{(j)}(X_{g1}) |\mathbf{X}\bigr\}\nonumber\\
&&\qquad= 2\sigma^2\{1-2\mathbf e_g^T(\mathbf I-\mathbf{S}_1^*\mathbf{S}_j^*)^{-1}
(\mathbf I-\mathbf{S}_1^*)\mathbf e_g\\
&&\hspace*{20.4pt}\qquad\quad{} +\mathbf e_g^T(\mathbf I-\mathbf{S}_1^*\mathbf{S}_j^*)^{-1}
(\mathbf I-\mathbf{S}_1^*)(\mathbf I-\mathbf{S}_1^*)^T(\mathbf I-\mathbf
{S}_1^*\mathbf{S}_j^*)^{-T}\mathbf e_g\}.\nonumber
\end{eqnarray}
Using the same argument as
that for (10) in Opsomer and Ruppert (\citeyear{OR97}), we obtain that
for $j,k=2,\ldots,J$ $(j\ne k)$,
\[
\operatorname{Cov}\bigl(\hat{m}_1^{(j)}(X_{g1}),\hat
{m}_1^{(k)}(X_{g1}) |\mathbf{X}\bigr)
=\frac{1}{Gh_1}\sigma^2f_1^{-1}(X_{g1})\nu_0(K)
+o_p\biggl(\frac{1}{Gh_1}\biggr)
\]
and
\[
\operatorname{Var}\bigl(\hat{m}_1^{(j)}(X_{g1}) |\mathbf{X}\bigr)
=\frac{2}{Gh_1}\sigma^2f_1^{-1}(X_{g1})\nu_0(K)
+o_p\biggl(\frac{1}{Gh_1}\biggr).
\]
Therefore, by (\ref{B5}),
\begin{eqnarray*}
\operatorname{Var}(\hat{m}_1(X_{g1}) |\mathbf{X})
&=&(J-1)^{-2}\sum_{j,k=2}^J \operatorname{Cov}\bigl(\hat
{m}_1^{(j)}(X_{g1}),\hat
{m}_1^{(k)}(X_{g1}) |\mathbf{X}\bigr)\\
&=&\frac{J}{J-1}\frac{1}{Gh_1}\sigma^2f_1^{-1}(X_{g1})\nu_0(K)+o_p\biggl(\frac
{1}{Gh_1}\biggr).
\end{eqnarray*}
The conditional bias of $\hat{m}_1$ is obviously the sum of biases for
each $\hat{m}_1^{(k)}$ ($k=2,\ldots,J$).
This completes the proof of the theorem.
\end{pf*}
\begin{pf*}{Proof of Theorem \protect\ref{thbb}}
The result can be proved along the line of Theorem~3.1 in Opsomer
(\citeyear{O00}).
\end{pf*}
\end{appendix}

\section*{Acknowledgments}
The authors thanks the Associate Editor and the referees for
constructive comments that substantially improved an earlier version of
this paper.

\printaddresses

\end{document}